\documentclass{article}
\usepackage[utf8]{inputenc}

\title{The role of inter-regional mobility \\ in forecasting SARS-CoV-2 transmission} 
\author{Martijn H. H. Schoot Uiterkamp \and Martijn G\"osgens \and Hans Heesterbeek \and Remco van der Hofstad \and Nelly Litvak}

\usepackage[numbers]{natbib}
\usepackage{graphicx}
\usepackage{subcaption}

\usepackage{amsmath}
\usepackage{amsfonts}
\usepackage{amssymb}
\usepackage{amsthm}
\usepackage[a4paper,margin=2.54cm]{geometry}
\usepackage{booktabs}
\usepackage[textsize=small]{todonotes}
\usepackage{comment}
\usepackage{tabularx}
\usepackage{hyperref}

\begin{document}

\maketitle

\begin{abstract}
    In this paper, we present a method to forecast the spread of SARS-CoV-2 across regions with a focus on the role of mobility. Mobility has previously been shown to play a significant role in the spread of the virus, particularly between regions. Here, we investigate under which epidemiological circumstances incorporating mobility into transmission models yields improvements in the accuracy of forecasting, where we take the situation in the Netherlands during and after the first wave of transmission in 2020 as a case study. We assess the quality of forecasting on the detailed level of municipalities, instead of on a nation-wide level. To model transmissions, we use a simple mobility-enhanced SEIR compartmental model with subpopulations corresponding to the Dutch municipalities. We use commuter information to quantify mobility, and develop a method based on maximum likelihood estimation to determine the other relevant parameters. We show that taking inter-regional mobility into account generally leads to an improvement in forecast quality. However, at times when policies are in place that aim to reduce contacts or travel, this improvement is very small.

\end{abstract}


\section{Introduction}

In the effort to contain a pandemic, local contact-restricting measures tailored to particular regions within a country can be beneficial when those regions exhibit large differences in infection prevalence. For instance, local measures can be implemented specifically for regions with a high level of infections, whereas such measures are not necessary for other regions with lower levels of infections. To take such decisions on a regional level, the local policy-makers need information on the effect that measures have on their specific region. Therefore, for such a region-focused approach to be successful, it is important to know the effect of measures on the progression of infections on this local level, rather than only on the national level.

Since the start of the COVID-19 pandemic, a substantial number of studies have focused on forecasting the course of the pandemic, for example to be able to judge the effect of policy changes (see, for example, the overview in \cite{Ray2020}). Already at the very start of the pandemic, the role of human mobility has been studied extensively, establishing a correlation between reductions in mobility (e.g., as a consequence of national restrictive measures) and the spread of the virus (see, e.g., \cite{Zhou2020,Gosgens2021,Wellenius2021}). A general conclusion from these studies is that human mobility between regions contributes significantly to the initial transmission and spread of the virus. Similarly, a reduction in mobility as a consequence of restrictive measures is shown to correlate with a reduction of new infections. This suggests that information on the influence of inter-regional mobility on the virus spread is useful for regional policy-makers.

When focusing on the influence of mobility on the  spread of the virus, it is important to identify which type of mobility is (most) relevant and which level of granularity and detail in mobility information is required. On the one hand, the country-level may be sufficient for inference and forecasting the general development of the pandemic, but is not suitable to gain insight in the emergence of local outbreaks or super-spreading events, which have been shown to be important \citep{Lemieux2021,Liu2020}. On the other hand, small-scale models, for example on the level of individuals, can provide enormous insight in demographic aspects of the transmission process \cite{Chang2021}, but require many assumptions on the behavior of individuals and such assumptions are generally hard to validate \cite{Fagiolo2007,Lorig2021}. Another reason why determining the right level of detail is important, is that exact mobility information is not available for reasons of privacy and thus proxies must be used instead. The availability and accuracy of such proxies depends on the level of granularity and detail, where smaller-scale and more detailed proxies are generally less readily available and more sensitive to noise. Summarizing, it is important to determine which level of granularity is sufficient to incorporate mobility in transmission models and which level of detail of mobility proxies is sufficient for this.

 A perhaps even more fundamental question is under which epidemiological conditions and restrictive measures does incorporating mobility in transmission models actually yield significant improvements in estimating the state of the epidemic and forecasting its development. For instance, several studies \citep{Badr2020,Gatalo2020,Nouvellet2021} suggest that the aforementioned correlation between mobility reduction and transmission is significantly weaker for time periods after the first wave as a consequence of different and changing sets of restrictive measures.

The goal of this paper is develop methodological tools to ultimately support regional policy-makers in taking appropriate contact-restrictive measures to contain the spread of SARS-CoV-2, where we take the epidemic in the Netherlands as a use-case. For this, we study the influence of inter-regional mobility on the spread of infections and aim to answer the following two questions:

\begin{enumerate}
    \item The first question is how effective particular (types of) measures are in reducing the number of local and inter-regional contacts . To answer this question, we estimate these numbers of contacts for our case study of the situation in the Netherlands using a mobility-enhanced SEIR-compartmental model \cite{Gosgens2021}. We find that measures that restricted visiting public places in August 2020 and a set of measures representing a partial lockdown in October and November 2020 (including the closing of bars and restaurants and a strong advice to work from home) were followed by a reduction in both the number of local and inter-regional contact moments. However, further restricting contacts by additional measures, such as minimizing the number of allowed guests at home and the closing of public places in November 2020, was followed by an increase in local contacts.
    
    \item The second question is whether and when taking inter-regional mobility into account improves the quality of short-term forecasts of the spread of infection. We answer this question by comparing the forecast quality of the aforementioned model to that of the same model without mobility. We find that taking mobility into account generally improves the forecast quality of the model. However, under measures that aim to reduce the amount of commuter travel and work contacts, this improvement turns out to be insignificant. This suggests that, under such measures, information on this type of travel is not necessary to explain the development of new infections.
\end{enumerate}

We build in this work on the mobility-enhanced SEIR-compartmental model of \cite{Gosgens2021}, where a given area is divided into several smaller regions (e.g., a country that is divided into municipalities). Transmissions within a given region are modeled using a standard SEIR-compartmental model \cite{Diekmann2012} and transmissions caused by contacts with inhabitants from other regions are modeled by enhancing the standard formulation of the SEIR-model with an additional term involving the mobility between these regions and their susceptible and infectious populations. 

Several parameters of this model can be taken from the literature since they concern the intrinsic properties of the virus that are not location-dependent (e.g., the average lengths of the latent and infectious periods). However, the transmission rates are dependent on the number of contact opportunities and are time-dependent. In earlier work \cite{Gosgens2021}, these parameter values were chosen based on the literature available at that time with the goal to serve as illustrative examples when assessing the trade-off between mobility restrictions and virus transmission. In several studies in the literature, a maximum likelihood (MLE) approach is followed to estimate the single transmission rate for compartmental models that do not include mobility \cite{Kirkeby2017}. In the current work, we develop a MLE procedure to determine both transmission rates in the mobility-enhanced model of \cite{Gosgens2021} simultaneously based on the infection prevalence of SARS-CoV-2 in the Dutch population. In order to do this, we turn the corresponding difference equations for the compartmental model used in \cite{Gosgens2021} into a partially {\em stochastic} model, where one of the difference equations is interpreted as the parameter of a suitable Poisson random variable.

Our evaluations indicate that the mobility-enhanced model and estimation procedure yield accurate forecasts of the distribution of infections throughout the country, i.e., the percentage of all reported infections that occurs in each region. In particular, we obtain accurate forecasts of the \emph{order} of regions in terms of their share in the total number of reported infections. This suggests that the model and the mobility data are sophisticated enough to describe the inter-regional patterns in virus transmission. This is important information for regional policy-makers since they are interested in the potential sources of local COVID-19 outbreaks. However, we also find that the model is able to accurately forecast the \emph{volume} of infections only during longer periods without changes in mobility-restricting measures. The predictive quality of the model reduces significantly within the first week after such changes but, in the subsequent week, returns to the original degree of accuracy. The main cause of this reduction in accuracy is the absence of information on the most recently infected individuals.

In this paper, we focus primarily on spread in the Netherlands. However, we believe that our general conclusions with regard to, for example, the added value of accounting for mobility also apply to similar countries, i.e., densely populated countries with substantial inter-regional mobility. Furthermore, although we focus here on the SARS-CoV-2 pandemic, our approach is also applicable to otherinfectious diseases with similar epidemiological characteristics.

The organization of the paper is as follows. In Section~\ref{sec_compartment}, we introduce and explain the compartmental model and in Section~\ref{sec_data}, we explain how we initialize the parameters of this model and which data sources we use for this. In Section~\ref{sec_methods}, we describe how we estimate the transmission rates using MLE, how we validate the model, and how we use it to investigate the influence of mobility. In Section~\ref{sec_results}, we present the results and, finally, in Section~\ref{sec_conclusion} we discuss the limitations of our work and state our conclusions.

\section{The compartmental model}
\label{sec_compartment}

In this section, we describe the compartmental model introduced in \cite{Gosgens2021} that we use to describe the spread of infection. We divide the population of size $N$ into a set $D$ of separate regions -- for example administrative regions such as municipalities or provinces -- and denote the population of each area $i \in D$ by $N_i$. Throughout, we assume that the population remains constant over time. At each time $t$, the population of each region $i$ is partitioned into six compartments that indicate their epidemiological state, denoted by $(S_i(t),E_i(t),I_i^T(t),I_i^U(t),R_i^T(t),R_i^U(t))$. These compartments contain all susceptible (not infected), exposed (infected but not yet infectious), positively tested infectious, untested infectious, recovered positively tested and recovered untested individuals, respectively. In our model, we distinguish between positively tested and untested individuals for two reasons. First, we assume that positively tested individuals do not travel to, or receive visitors from, other regions. Second, the compartment of positively tested infectious persons is the only compartment whose size we can accurately measure or estimate.

We model mobility between regions as follows. For each two regions $i,j \in D$, we let the parameter $M_{ij}$ denote the number of individuals traveling from $i$ to $j$ per time unit. We assume that traveling individuals visit only one region and return directly afterwards to their home region. This means, among others, that no traveling person is counted twice in the parameters $M_{ij}$ and that travelers have contacts only with individuals in their destination region. For convenience of notation, we assume that $M_{ii} = 0$ for all $i \in D$.

In compartmental epidemiological models, the transmission rate, usually denoted by $\beta$, determines the rate at which susceptible individuals are infected by infectious ones. This rate depends on the number of close contacts that people have that could potentially lead to transmission: the more such contacts by an infectious individual, the higher the expected number of transmissions. With the inclusion of mobility between regions, two separate transmission rates $\beta_{\text{loc}}$ and $\beta_{\text{mob}}$ are introduced that correspond to the transmission rates for infections caused by intra-regional contacts and inter-regional contacts, respectively.

We assume that infections among individuals from a given region $i$ can occur via three types of contacts:
\begin{enumerate}
    \item Contacts between a susceptible and an infectious individual both belonging to region $i$ (local transmissions);
    \item Contacts between a susceptible individual from region $i$ that is traveling to another region $j$ and an infectious (untested) individual from $j$ (mobility-induced transmission due to travellers meeting infectious individuals);
    \item Contacts between a susceptible individual from region $i$ that is visited by an infectious (untested) individual from another region $j$ (mobility-induced transmission due to infectious travellers meeting healthy individuals).
\end{enumerate}
This means that the total number of new infections among persons belonging to region $i$  at time $t$ caused by local contacts within $i$ is given by
\begin{equation*}
    \beta_{\text{loc}} \frac{S_i(t)}{N_i} (I_i^{\text{T}}(t) + I_i^{\text{U}}(t)),
\end{equation*}
where the ratio $S_i(t)/N_i$ can be interpreted as the probability that a random contact of an individual in region $i$ is with a susceptible individual. Here, we assume that the population within the region is \emph{well-mixed}, meaning that each pair of individuals has the same probability of coming into contact in the appropriate manner.

The number of new infections in region $i$ caused because susceptible individuals from $i$ travel to other regions is given by
\begin{equation*}
    \beta_{\text{mob}} \sum_{j \in D} M_{ij} \frac{S_i(t)}{N_i} \frac{I_j^{\text{U}}(t)}{N_j},
\end{equation*}
while the number of new infections caused because infectious individuals from other regions travel to region $i$ is given by
\begin{equation*}
    \beta_{\text{mob}} \sum_{j \in D} M_{ji} \frac{S_i(t)}{N_i} \frac{I_j^{\text{U}}(t)}{N_j}.
\end{equation*}
In both expressions, the ratio $(S_i(t)/N_i) (I_j^U(t)/N_j)$ can be interpreted as the probability that a contact between individuals from regions $i$ and $j$ involves a susceptible person from $i$ and an infectious person from $j$. Note that we assume that positively-tested infectious persons do not have contact with travelers and thus do not contribute to transmissions between individuals from different regions.

Summarizing, the compartmental model is described by the following differential equations, for each $i \in D$:
\begin{align}
\frac{dS_i(t)}{dt} &= -\beta_{\text{loc}} \frac{S_i(t)}{N_i} (I_i^{\text{T}}(t) + I_i^{\text{U}}(t))
-\beta_{\text{mob}} \sum_{j \in D} \left( S_i(t) \frac{M_{ji}}{N_i} \frac{I_j^{\text{U}}(t)}{N_j} + I_j^{\text{U}}(t) \frac{M_{ij}}{N_j} \frac{S_i(t)}{N_i} \right), \label{dfeq_01} \\
\frac{dE_i(t)}{dt} &= \beta_{\text{loc}} \frac{S_i(t)}{N_i} (I_i^{\text{T}}(t) + I_i^{\text{U}}(t))
+\beta_{\text{mob}} \sum_{j \in D} \left( S_i(t) \frac{M_{ji}}{N_i} \frac{I_j^{\text{U}}(t)}{N_j} + I_j^{\text{U}}(t) \frac{M_{ij}}{N_j} \frac{S_i(t)}{N_i} \right) - \frac{E_i(t)}{\nu}, \nonumber \\
\frac{dI_i^{\text{T}}(t)}{dt} &= a \frac{E_i(t)}{\nu} - \frac{I_i^{\text{T}}(t)}{\omega}, \label{dfeq_03} \\
\frac{dI_i^{\text{U}}(t)}{dt} &= (1-a) \frac{E_i(t)}{v} - \frac{I_i^{\text{U}}(t)}{\omega}, \nonumber \\
\frac{dR_i^{\text{T}}(t)}{dt} &= \frac{I_i^{\text{T}}(t)}{\omega}, \nonumber \\
\frac{dR_i^{\text{U}}(t)}{dt} &= \frac{I_i^{\text{U}}(t)}{\omega}. \nonumber
\end{align}
Here, $a$ is the fraction of infectious individuals that is positively tested, $\nu$ is the average latent period, i.e., the time between getting infected and becoming infectious, and $\omega$ is the average infectious period, during which an infected individual can infect others.

The transmission rates implicitly contain information on the number of contacts that people have and on the infectiousness of the virus. To analyze changes in these values over the course of the pandemic, we decompose the transmission rates in terms of these quantities. More precisely, observe that these rates can be interpreted as the product of the probability that a transmission occurs given an appropriate contact between a susceptible and an infectious person, denoted by $\varepsilon$, and the average number of such contacts per person per unit time within their own area (for $\beta_{\text{loc}}$) and when traveling (for $\beta_{\text{mob}}$). Given the overall average contact rate $c$ and the fraction $p$ of those contacts that occur locally, i.e., within a given area, we can decompose $\beta_{\text{loc}}$ and $\beta_{\text{mob}}$ as
\begin{align}
    \beta_{\text{loc}} &= \varepsilon p c, \label{eq_beta_loc}\\
    \beta_{\text{mob}} &=  \varepsilon (1-p) \frac{cN}{2 \sum_{(i,j) \in D \times D} M_{ij}}. \label{eq_beta_mob}
\end{align}
In the latter expression, the term 
    \[
    (1-p) \frac{cN}{2 \sum_{i,j \in D \times D} M_{ij}}
    \]
represents the total number of non-local contacts $(1-p)\frac{cN}{2}$ divided by the total number of travels $\sum_{(i,j) \in D^2} M_{ij}$, i.e., the number of contacts per traveling person per unit time.

A concise explanation of the model parameters is given in Table~\ref{tab_par}. The parameters $D$ and $N_i$ are initialized based on the division of the country into regions. For instance, if we investigate on the municipality-level, then $D$ contains an index for each municipality. Furthermore, the values for the average latent and infectious periods $\nu$ and $\omega$ are chosen based on the literature as 3 days \citep{He2020,Li2020} and 9 days \citep{Byrnee2020,Park2020,Walsh2020}, respectively. We assume that all above parameters are constant over time. In the following section, we discuss how the remaining time-dependent parameters are initialized. 

\begin{table}[ht!]
\centering
\begin{tabularx}{\textwidth}{l X}
\toprule
Parameter & Meaning \\
\midrule
$N$ & Total population size \\
$D$ & Set of region indices  \\
$N_i$ & Population size of region $i$ \\
$M_{ij}$ & Number of individuals that travel from region $i$ to region $j$ \\
$a$ & Fraction of infectious people that have been tested \\
$\omega$ & Infectious period \\
$v$ & Latent period \\
$\beta_{\text{loc}}$ & Transmission rate via local contacts \\
$\beta_{\text{mob}}$& Transmission rate via non-local contacts \\
$\varepsilon $ & Transmission probability per contact between susceptible and infectious individual \\
$c$ & Average number of contacts \\
$p$ & Fraction of local contacts \\
\bottomrule
\end{tabularx}
\caption{Model parameters.}
\label{tab_par}
\end{table}

\section{Parameter initialization and data sources}
\label{sec_data}
 In this section, we describe how we initialize parameters in the compartmental model. Respectively, we focus on the initialization of the compartments themselves, the mobility parameters $M_{ij}$, the fraction of positively tested infectious individuals $a$, the fraction of local contacts $p$, and the product $\varepsilon c$ of the infection probability and the number of contacts. Note that we do not focus on initializing the parameters $\varepsilon$ and $c$ separately since they do not occur separately in our model but only together in the form $\varepsilon c$.

\subsection{Initialization of compartments}
\label{sec_compartments}

To initialize the sizes of the six compartments on a given day $t$ for a given municipality $i$, we use information on daily reported cases as reported by, for example official national health sources such as the Dutch National Institute for Public Health and the Environment (RIVM) \cite{RIVM}. More precisely, we use as input the number of positive tests reported at time $t$ in region $i$, denoted by $\Delta I_i(t)$, for each time $t$ and region $i$. The corresponding individuals have been infected before $t$ and are either infectious, not yet infectious, or not infectious anymore. However, we assume that recipients of a positive test result become infectious at the report date of the test since we do not have (access to) any other information on when a positively tested individual has been infected, underwent a test, or recovered.

By the structure of the compartmental model, the latent period of an exposed individual has an exponential distribution with expected value $\nu$. This implies for any $s \geq 0$ that the expected number of positively tested individuals that were exposed at time $t$ that become infectious at time $t+s$ is $e^{-\frac{s}{\nu}} \frac{\Delta I_i(t+s)}{a}$. We initialize the number of exposed individuals as at time $t$ as
\begin{equation}
    \hat{E}_i(t) := \frac{1}{a} \int_0^{\infty}e^{-\frac{s}{\nu}}
    \Delta I_i(t+s) ds.
    \label{eq_E}
\end{equation}
Note that this approach requires information on \emph{future} positive tests.

Analogous to the latent period, the infectious period of an individual has an exponential distribution with expected value $\omega$. This means for any $s \geq 0$ that the expected number of tested individuals that become infectious at time $t-s$ and are still infectious at time $t$ is $e^{-\frac{s}{\omega}} \Delta I_i(t-s)$. We initialize the number of tested infectious and recovered individuals at time $t$ as
\begin{equation}
    \hat{I}_i^{\text{T}}(t) :=
    \int_0^{\infty}e^{-\frac{s}{\omega}}
    \Delta I_i(t-s) ds,
    \label{eq_It}
\end{equation}
and
\begin{equation}
    \hat{R}_i^{\text{T}}(t) :=
    \int_0^{\infty} 
    \Delta I_i(t-s) ds
    - \hat{I}_i^{\text{T}}(t),
    \label{eq_Rt}
\end{equation}
respectively.

Given the fraction $a$ of tested infectious individuals, we initialize the compartments $I_i^{\text{U}}(t)$ and $R_i^{\text{U}}(t)$ for the number of untested infectious and recovered persons by scaling the compartments $I_i^{\text{T}}(t)$ and $R_i^{\text{T}}(t)$ by the factor $\frac{1-a}{a}$ (we explain how we estimate $a$ in Section~\ref{sec_frac_tested}). To see this, observe that $I_i^{\text{U}}(t) = (1-a)(I_i^{\text{T}}(t) + I_i^{\text{U}}(t))$, which implies that $\hat{I}_i^{\text{U}}(t) = \frac{1-a}{a} \hat{I}_i^{\text{T}}(t)$. Similarly, it follows that $\hat{R}_i^{\text{U}}(t) := \frac{1-a}{a} \hat{R}_i^{\text{T}}(t)$. Finally, the number of susceptible persons $S_i(t)$ follows directly from the number of inhabitants $N_i$ and the sizes of the other compartments as
\begin{equation*}
    \hat{S}_i(t) := N_i - \hat{E}_i(t) - \hat{I}_i^{\text{T}}(t) - \hat{I}_i^{\text{U}}(t) - \hat{R}_i^{\text{T}}(t) - \hat{R}_i^{\text{U}}(t).
\end{equation*}

\subsection{Initialization of mobility}
\label{sec_mobility}

We now describe how the mobility parameters $M_{ij}$ are initialized. For this initialization, we combine two different data sources: one on mobility patterns from before the pandemic and the other on the relative change in mobility since the start of the pandemic.

The first source is data on places of residence and of work as collected by the Statistics Netherlands (CBS) \cite{CBS2019}. For each two Dutch municipalities, this data specifies the number of employees that live in the former municipality and work at a company registered in the latter municipality. Thus, these data provide a proxy for the structure of commuting travel in the Netherlands. The most recently collected version of these data, which we also use in this work, is from December 2019 and contains commuting information between all 355 Dutch municipalities.

Clearly, both the volume and structure of mobility patterns has changed significantly since the start of the pandemic \citep{deHaas2020}. In particular, both the total number of travel movements and the share of long-distance trips among the total number of movements have significantly decreased \citep{Schlosser2020}. As a consequence, simply using the mobility proxy obtained from the CBS data as mentioned above may not be sufficient to make accurate estimations and forecasts. 

To account for changes in mobility that have occurred as a consequence of the pandemic, we reduce the commuting behavior as given by the CBS data by using Google mobility data \citep{Google2021}. These data represent the change in mobility for different types of mobility compared to a baseline established in January 2020. More precisely, a division is made between mobility involving different destination categories such as retail, residential, and workplaces. The baseline consists of mobility patterns on seven days, one for each day of the week. For each of these days, the baseline day is constructed by taking the median mobility of these days within the 5-week period January 3 - February 6, 2020.

To obtain a proxy $\hat{M}$ for the mobility at a given moment after March 2020, we use the CBS data as a proxy for the baseline of the Google mobility data. This means that we directly scale the CBS commuting data by the rates given in the Google data with regard to workplace mobility. For the Netherlands, the Google mobility data contains reduction rates on three spatial granularity levels, namely on the level of municipalities, provinces, and the country as a whole. Missing entries on the municipality level were substituted by the corresponding entry of the province to which the given municipality belongs. Analogously, missing entries on the province level were substituted by the corresponding entry on the country level.

\subsection{Initializing the fraction of positively tested infectious people}
\label{sec_frac_tested}

Via studies on the seroprevalence of SARS-CoV-2 in the population, the number of hospitalizations, and the daily fraction of positive tests, estimates can be made of the total number of infectious people $I^{\text{tot}}(t):= \sum_{j \in D} (I_j^{\text{T}}(t) + I_j^{\text{U}}(t))$ at time $t$. Based on such an estimate $\hat{I}^{\text{tot}}$ and the total number of positively tested individuals as initialized in Section~\ref{sec_compartments}, we can obtain an estimate on the fraction of positively tested infectious individuals $a$ on a given day $t$ as
\begin{equation*}
    a = \frac{\sum_{j \in D} I_j^{\text{T}}(t)
    }
    {\sum_{j \in D} (I_j^{\text{T}}(t) + I_j^{\text{U}}(t))
    }
    =
    \frac{ \sum_{i \in D}
     \int_0^{\infty}e^{-\frac{s}{\omega}}
    \Delta I_i(t-s) ds
    }
    {\hat{I}^{\text{tot}}(t)}.
\end{equation*}

For $I^{\text{tot}}(t)$ we use estimates  by the RIVM based on the outcomes of the ``Pienter Corona-studie'' \citep{Vos2020}, which investigated the seroprevalence of SARS-CoV-2 in the Netherlands (see also \cite{Verberk2019}).

\subsection{Initializing the fraction of local contacts $p$ and infection probability times number of contacts $\varepsilon c$}
Given the transmission rates $\beta_{\text{loc}}$ and $\beta_{\text{mob}}$, we estimate the fraction of local contacts $p$ as follows (we describe how we estimate the transmission rates in Section~\ref{sec_MLE}). Note that, by the definition of $\beta_{\text{loc}}$ and $\beta_{\text{mob}}$ in Equations~(\ref{eq_beta_loc}) and~(\ref{eq_beta_mob}),
\begin{equation*}
\frac{\beta_{\text{loc}}}{p}
= \varepsilon c
= \frac{2 \sum_{(i,j) \in D \times D} M_{ij}}{(1-p) N} \beta_{\text{mob}},
\end{equation*}
provided that $p \neq 0$ and $p \neq 1$. From this, it follows that
\begin{equation*}
\frac{1}{p} = \frac{2 \sum_{(i,j) \in D \times D} M_{ij}}{N} \frac{\beta_{\text{mob}}}{\beta_{\text{loc}}} + 1,
\end{equation*}
and thus
\begin{equation*}
    p = \frac{N \beta_{\text{loc}}}{2 \beta_{\text{mob}} \sum_{(i,j) \in D \times D} M_{ij} + N \beta_{\text{loc}}}.
\end{equation*}
Observe that we can now compute the term $\varepsilon c$ as $\varepsilon c= \beta_{\text{loc}}/p$.

\section{Approach}
\label{sec_methods}

In this section, we describe how we use the compartmental model and the data initialization procedure in order to accomplish the following:
\begin{enumerate}
    \item estimate the transmission rates $\beta_{\text{loc}}$ and $\beta_{\text{mob}}$;
    \item validate the performance of the model at the nationally-aggregated level;
    \item assess the forecasting quality of the model on the municipality level;
    \item investigate the influence of mobility on the forecast quality.
\end{enumerate}
Before this, we briefly discuss how we discretize the model so that it matches the discrete-time format of the available input data.

We focus on the developments within the second half of 2020, i.e., within the time period July 1 - December 31 2020. The reason for this is that, on the one hand, testing for the virus became available for the entire population slightly before July 1, 2020. Thus, from this time on, the daily reported cases provide a more accurate view of the number of known infections. On the other hand, mass vaccination of the population started in January 2021. Thus, the current version of our compartmental model is especially suitable for the time period July 1, 2020 - December 31, 2020.

\subsection{Model discretization}
The data that is used to initialize the compartments is reported on a discrete (daily) timescale. Therefore, we discretize the time domain in days and simulate the compartmental model via an appropriate finite differences method. To initialize the compartments, we discretize the integrals in the expressions for $\hat{E}_i(t)$, $\hat{I}_i^{\text{T}}(t)$, and $\hat{R}_i^{\text{T}}(t)$ in Equations~(\ref{eq_E}-\ref{eq_Rt}). We replace the exponential distributions for the lengths of the latent and infectious periods by geometric distributions with success probabilities $\frac{1}{\nu}$ and $\frac{1}{\omega}$, respectively, that capture the same memory-less behavior. These parameters are chosen such that the expectation remain $\nu$ and $\omega$, respectively. Moreover, we truncate the resulting sums so that we look ahead 7 days (when calculating $\hat{E}_i(t)$) and look back 14 days (when calculating $\hat{I}_i^{\text{T}}(t)$ and $\hat{R}_i^{\text{T}}(t)$), where we put the remainder probability mass on the last day within this look-ahead or look-back time window. This means that we calculate $\hat{E}_i(t)$, $\hat{I}_i^{\text{T}}(t)$, and $\hat{R}_i^{\text{T}}(t)$ as
\begin{align*}
    \hat{E}_i(t) &= \sum_{s=1}^{6} \left(1 - \frac{1}{\nu} \right)^{s-1} \Delta I_i(t+s) +  \left(\nu  - \nu \left(1 - \frac{1}{\nu}\right)^6 \right) \Delta I_i(t + 7); \\
    \hat{I}_i^{\text{T}}(t) &= \sum_{s=0}^{13} \left(1-\frac{1}{\omega} \right)^s \Delta I_i(t-s) + \left(\omega - \omega \left(1 - \frac{1}{\omega} \right)^{14} \right) \Delta I_i(t-14); \\
    \hat{R}_i^{\text{T}}(t) &= \sum_{s=0}^{t - s_{\text{start}}} \Delta I_i(t-s) - \hat{I}_i^{\text{T}}(t),
\end{align*}
where $s_{\text{start}}$ is the first day in the pandemic on which infections were reported.

\subsection{Estimating the transmission rates $\beta_{\text{loc}}$ and $\beta_{\text{mob}}$}
\label{sec_MLE}

 We develop an MLE procedure to estimate for a given day $t$ the corresponding transmission rates. This means that we assume that for a given day $t$ and municipality $i$ the number of new infections on the next day, i.e.,  $-\Delta S_i(t) := S_i(t) - S_i(t+1)$, is a random variable with given distribution. The goal is to find the parameter values of this distribution that make this assumption most likely.

In particular, in line with other work on MLE for estimating transmission rates in compartmental models (see, e.g., \cite{post2020,O'Driscoll2021}), we investigate two different but commonly used models for the distributions of $-\Delta S_i(t)$. In the first model, we assume that $-\Delta S_i(t)$ is a Poisson random variable whose rate is given by the negative of the right-hand side of Equation~(\ref{dfeq_01}), i.e., by
\begin{equation}
\lambda_i (t) := \beta_{\text{loc}} \frac{S_i(t)}{N_i} (I_i^{\text{T}}(t) + I_i^{\text{U}}(t))
+\beta_{\text{mob}} \sum_{j \in D} \frac{S_i(t)}{N_i} \frac{I_j^{\text{U}} (t)}{N_j} (M_{ji} + M_{ij}) .
\label{eq_MLE_rate}
\end{equation}

In the second model, we assume that $-\Delta S_i(t)$ has a negative binomial distribution whose mean is given by $\lambda_i(t)$ as defined in Equation~(\ref{eq_MLE_rate}) and whose variance is given by $\lambda_i(t) (1 + \frac{\lambda_i(t)}{r(t)})^{r(t)}$, where $r(t)$ is a to-be-estimated dispersion parameter. In contrast to the Poisson distribution, the negative binomial distribution models overdispersion since it can be interpreted as a mixed Poisson random variable whose parameter has a Gamma distribution. Therefore, it offers a suitable alternative to the Poisson model when the underlying data exhibits higher variability than a Poisson variable.

Note that in both models, only the flow from $S_i(t)$ to $E_i(t)$ is taken as a stochastic variable and not the flows between other compartments as is common in some other works (see, e.g., \cite{ShaHe2020}). The main reason for this is that these flows depend on parameters that are necessary for the compartment initialization and must therefore be estimated beforehand in a different way. Thereby, the total population size remains constant for each municipality.

The maximum likelihood approach relies on estimating parameters by optimizing an appropriate (log-)likelihood function, which we describe next. For the case of Poisson random variables, the log-likelihood function for the parameters $\beta_{\text{loc}}$ and $\beta_{\text{mob}}$ at time $t$ is given by
\begin{align*}
& L_{\text{P}} \left( 
\beta_{\text{loc}},\beta_{\text{mob}} \ | \ x_{\text{loc}}(t),x_{\text{mob}}(t) , \Delta S(t)
\right)\\
:= &
\sum_{i \in D} \left(
-\Delta S_i(t) \log (\beta_{\text{loc}} x_{i,\text{loc}}(t)
+ \beta_{\text{mob}} x_{i,\text{mob}}(t))
- \beta_{\text{loc}} x_{i,\text{loc}}(t)
- \beta_{\text{mob}} x_{i,\text{mob}}(t)
-
\log((-\Delta S_i(t)) !)
\right),
\end{align*}
where the $i^{\text{th}}$ elements of the vectors $x_{\text{loc}}(t)$ and $x_{\text{mob}}(t)$ are given by
\begin{align*}
x_{i,\text{loc}}(t) & := \frac{S_i(t)}{N_i}(I_i^{\text{T}}(t) + I_i^{\text{U}}(t)), \\
x_{i,\text{mob}}(t) & := \sum_{j \in D} \frac{S_i(t)}{N_i} \frac{I_j^{\text{U}} (t)}{N_j} (M_{ji} + M_{ij}) .
\end{align*}
For the case of negative binomial variables, the log-likelihood function is given by
\begin{align}
    &L_{\text{NB}} \left(\beta_{\text{loc}},\beta_{\text{mob}},r \ | \ x_{\text{loc}}(t), x_{\text{mob}}(t),\Delta S(t) \right) \nonumber \\
    =& \sum_{i \in D} \log \left(
    \frac{(\beta_{\text{loc}} x_{i,\text{loc}}(t)
+ \beta_{\text{mob}} x_{i,\text{mob}}(t))^{-\Delta S_i(t)}}{(-\Delta S_i(t))!}
    \frac{\Gamma(r - \Delta S_i(t))}{\Gamma(r) (r + \beta_{\text{loc}} x_{i,\text{loc}}(t)
+ \beta_{\text{mob}} x_{i,\text{mob}}(t))^{-\Delta S_i(t)}} \right. \nonumber \\
& \quad \cdot \left.
    \frac{1}{(1 + \frac{\beta_{\text{loc}} x_{i,\text{loc}}(t)
+ \beta_{\text{mob}} x_{i,\text{mob}}(t)}{r})^r}
    \right) \nonumber \\
    =&
    \sum_{i \in D} \left( -\Delta S_i \log(\beta_{\text{loc}} x_{i,\text{loc}}(t)
+ \beta_{\text{mob}} x_{i,\text{mob}}(t)) - \log((-\Delta S_i(t))!) 
    + \log(\Gamma(r - \Delta S_i(t)))
    - \log(\Gamma(r)) \right. \nonumber \\
    & \quad \left.
    + \Delta S_i(t) \log(r + \beta_{\text{loc}} x_{i,\text{loc}}(t)
+ \beta_{\text{mob}} x_{i,\text{mob}}(t))
    - r \log \left(1 + \frac{\beta_{\text{loc}} x_{i,\text{loc}}(t)
+ \beta_{\text{mob}} x_{i,\text{mob}}(t)}{r} \right) \right),
\label{eq_lll}
\end{align}
where $\Gamma$ is the Gamma-function. The desired estimate of $(\beta_{\text{loc}},\beta_{\text{mob}})$, and an estimation of the corresponding $r$ in the case of $L_{\text{NB}}$, is a maximizer of this function wherein we substitute our estimates for the compartments and mobility data as described in Section~\ref{sec_compartments}. Since both $L_{\text{P}}$ and $L_{\text{NB}}$ are concave in $\beta_{\text{loc}}$ and $\beta_{\text{mob}}$ (and $r$), standard descent algorithms can be employed to find such a maximizer efficiently.

Note that one consequence of estimating the transmission rates in this way is that the terms $\hat{\beta}_{\text{mob}} x_{i,\text{mob}}$ are invariant to scaling $M$ by a positive scalar. Thus, when the used proxy $\hat{M}$ for $M$ has the same structure as $M$ but differs primarily in the volume, i.e., if $\hat{M} \approx sM$ for some positive scalar $s$, the obtained estimated transmission rates are also suitable estimates for the transmission rates obtained when having used the true mobility matrix $M$. This also holds for the resulting estimate of the fraction of local contacts $p$. Note that this observation is independent of the data source used for the proxy $\hat{M}$. This observation explains why our methods are not very sensitive to changes in mobility volumes.

The confidence intervals for the estimates are computed via parametric bootstrapping. More precisely, we use the initial point estimates of the transmission rates and dispersion parameter to draw 100 random samples of the increase $-\Delta S_i$ in new infections from the negative binomial distribution whose rate is given by $\lambda_i$ in Equation~(\ref{eq_MLE_rate}) after substituting the point estimates for the transmission rates. Subsequently, we apply MLE for each of these samples to obtain 100 point estimates of the transmission rates and dispersion parameter that we use to construct the confidence intervals. Using these point estimates, we also construct confidence intervals for all model parameters that are derived from the transmission rates (the fraction of local contacts $p$ and the term $\varepsilon c$). Moreover, when forecasting new infections, we directly obtain 100 forecast scenarios from these estimates.

In initial testing, we found that the negative binomial model outperforms the Poisson model on practically all outputs. Therefore, in the remainder, we focus only on the negative binomial model.

\subsection{Validation on the national level}
\label{sec_val_nat}

The model validation consists of two steps. In each step, we compare an output of the model to the corresponding reported values by RIVM \cite{RIVM}. The two considered outputs are the daily reported positive tests and the effective reproduction number.

In the first validation step, we compare the daily reported positive cases on the national level as computed by our model with those reported by the RIVM. Within our model, the number of daily reported positive tests for a given day $t$ is represented via Equation~(\ref{dfeq_03}) by the term $\sum_{i \in D} a \frac{\hat{E}_i(t)}{\nu}$. We focus on two methods for calculating these terms from our model:
\begin{enumerate}
    \item
The first method is based purely on the initialization approach. More precisely, we directly compute for each day $t$ the corresponding term using Equation~(\ref{eq_E}) as
\begin{equation*}
    \sum_{i \in D} a \frac{\hat{E}_i(t)}{\nu}
    = \frac{1}{\nu} \int_0^{\infty}e^{-\frac{s}{\nu}}
    \Delta I_i(t+s) ds.
\end{equation*}
This method allows us to check, among others, the validity of the assumption that the latent and infectious periods have exponential distributions. 
\item
In the second method, for each day we first estimate the transmission rates as described in Section~\ref{sec_MLE}. Subsequently, we initialize the compartments at the first date, July 1, and simulate the model until the last day, December 31, using for each day the corresponding estimated transmission rates, fraction of positive tests, and mobility data as input. This procedure directly yields the desired terms for each simulated date. To reduce the influence of propagating errors, we re-initialize the compartments at every first day of the month. This method allows us to check the validity of the model dynamics, i.e., whether the model is able to adequately capture the development of new positive cases over time.
\end{enumerate}

In the second validation step, we use our model to estimate the effective reproduction number and compare these estimates to those reported by RIVM (see \cite{Wallinga2007,vandeKasteele2019} for a description of the estimation methods used by RIVM). For each day, we compute the effective reproduction number by constructing the \emph{next-generation-matrix} (NGM) \citep{Diekmann2010} based on our model and estimates of the transmission rates (for more details on the construction of this matrix, we refer to \cite{Gosgens2021}). If the general trend of our and RIVM's estimate matches, then the underlying methods to obtain these estimates agree on the representation of the dynamics between susceptible and infectious individuals. 

However, we also expect that, even if the general trends of our estimates and those of RIVM match, our estimates may have a delay of several days for two reasons. First, in our model we assume that tested individuals receive the result of their test immediately after they become infectious. More precisely, we do not account for the time between getting tested (for example because one develops symptoms) and receiving a (positive) test result. However, throughout the studied time period July 1 - December 31, the average time between planning a test appointment and receiving the result has varied between 35 and 81 hours in the Netherlands \citep{Rouvoet2021}. As a consequence, we might estimate the number of infectious and tested individuals with some delay, which is caused by the time between getting tested and receiving the result of a test. In contrast, the calculation method used by RIVM does allow for delays in  testing and reporting \citep{Wallinga2007,vandeKasteele2019}.

The second reason is that, to compute the effective reproduction number for a given date, our computation via the NGM and estimated transmission rates uses information on daily reported positive cases from both before and after the date as input. On the other hand, the aforementioned calculation method of RIVM utilizes detailed information about new reported hospital admissions. These are reported by hospitals with some delay and these numbers represent infections that took place up to 14 days previously. One reason that the RIVM uses hospital admissions is that the daily numbers of positive tests depend on the number of people that get tested and this number varies also for reasons other than variation in new cases (for example, changes in testing strategy in the early phases of the outbreak in the Netherlands.

To assess whether our estimates have some delay compared to the estimates of RIVM, we determine the amount of days by which we should shift our estimates as to obtain the largest correlation between these time-shifted estimates and those of RIVM.

\subsection{Forecast quality}
\label{sec_setup_forecast}

To assess the forecasting quality of our approach, we compute forecasts of the number of newly reported infections in the next 14 days. To obtain these forecasts, we initialize the compartmental model 7 days in the past and simulate it over a horizon of 21 days (recall that we cannot initialize the model for the current date since the number of exposed people is calculated from the number of reported cases in the next 7 days). We choose the necessary parameters $\beta_{\text{loc}}$, $\beta_{\text{mob}}$, and $a$ as the average of these parameters over the past 7 days from the start of the forecast horizon. As an example, suppose that the current date is August~11. Then we aim to make forecasts of the number of reported infections in the period August 11-25 and initialize the model on August~4 using as parameter inputs the averages of these parameters over the period July~28 - August~3.

We distinguish between the absolute and relative number of reported positive cases, where the relative number is in comparison to the total nation-wide number of infections. The latter values provide insight in the risk that infections will occur in a given municipality. Thus, they are useful risk information for local policy-makers. 

With respect to the absolute number of positive tests, we directly compare our forecasts to the actual number of reported cases within the chosen period, both on the national level and that of municipalities. With regard to the level of municipalities, we assess whether there is a spatial discrepancy in forecast accuracy between different municipalities. With respect to the relative number of positive tests, we investigate whether we can accurately forecast the \emph{order} of municipalities in terms of the number of reported positive tests. We do this by calculating the correlation between the \emph{orders} of the forecast and observed fractions of reported infections via the Spearman correlation.

\subsection{Influence of mobility}

We assess the influence of mobility in two ways. First, we investigate under which conditions including mobility leads to a better fit of the transmission rates in the negative binomial model. To this end, we consider two different models for estimating the transmission rates via MLE, which are given by the presence or absence of mobility. When estimating the reates when mobility is excluded, we set $\beta_{\text{mob}} = 0$ in the expressions for the log-likelihood function in Equation~(\ref{eq_lll}). To assess which of these models fits the data best, we compute and compare for both models the Akaike information criterion (AIC) \cite{Akaike1998}. For each model, this criterion is computed as $2k - 2\log \hat{L}$, where $k$ is the number of unknown parameters that the model estimates and $\log \hat{L}$ is the maximum value of the log-likelihood function. A lower AIC implies that the given model fits the data better. In particular, a difference in the AIC of two models of more than $10$ is generally considered to be a strong indication that the model with the lower AIC value is a better fit to the data \citep{Burnham2002}.

Second, we investigate under which conditions the inclusion of mobility leads to more accurate forecasts of the number of reported positive tests. For this, we compute for both model versions (including or excluding mobility) forecasts as described in the previous subsection and compare the difference in accuracy by means of the root mean squared error (RMSE) (for the absolute number of infections) and the Spearman correlation (for the relative number of infections).

\section{Results}
\label{sec_results}

We now present and discuss the results of the estimations and experiments described in the previous section\footnote{The corresponding source code is available at \url{https://github.com/MartijnGosgens/mobility-forecasting-covid}.}. To assess the influence of new measures on the estimated parameters and on the forecasting accuracy, we provide an overview of the restrictive measures taken in the Netherlands within the studied time period in Table~\ref{tab_timeline}. In all figures in this section, the dates at which these measures were installed are indicated by green (relaxation), orange (restriction), and red (lockdown) vertical lines.

\begin{table}[ht!]
\centering
\begin{tabularx}{\textwidth}{ll X}
\toprule
Date & Type & Decription \\
\midrule
July 1 & Relaxation & Only basic restrictions (keeping 1.5 meters distance, washing hands, etc.) \\
August 6 & Restriction & Stricter rules for attending public places (e.g., mandatory pre-registration) \\
August 18 & Restriction & At most 6 guests at home \\
September 29 & Restriction & At most 3 guests at home, bars and restaurants close at 22:00h, no audience at sport events; advice to work from home \\
October 14 & Partial lockdown & Bars and restaurants closed, no events, face masks mandatory at schools, sport in groups of at most 4; strong advice to work from home \\
November 4 & Partial lockdown & At most 2 guests at home, public places closed, sport in groups of at most 2 \\
December 14 & Lockdown & Schools, non-essential stores, and gyms closed \\
December 24-26 & Temporary relaxation & At most 3 guests at home \\
\bottomrule
\end{tabularx}
\caption{Overview of measures taken in the Netherlands between July 1 - December 31, 2020. Taken from \cite{tijdlijn}.}
\label{tab_timeline}
\end{table}

\subsection{Estimation of the transmission rates and related parameters}

\begin{figure}[ht!]
\centering
\includegraphics[scale = .7]{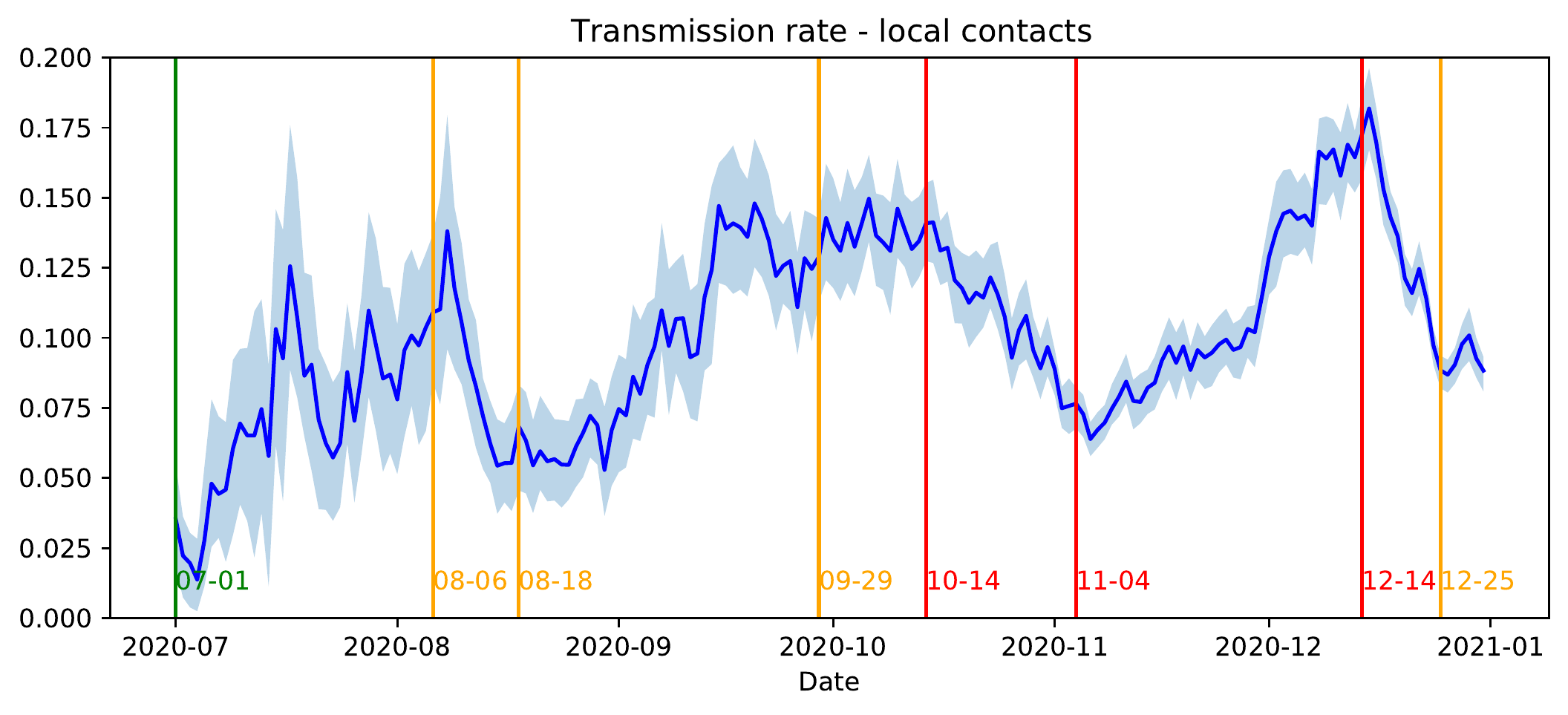}
\caption{Estimations of the local transmission rate $\beta_{\text{loc}}$ with 95\% confidence interval.}
\label{fig_loc}
\end{figure}

\begin{figure}[ht!]
\centering
\includegraphics[scale = .7]{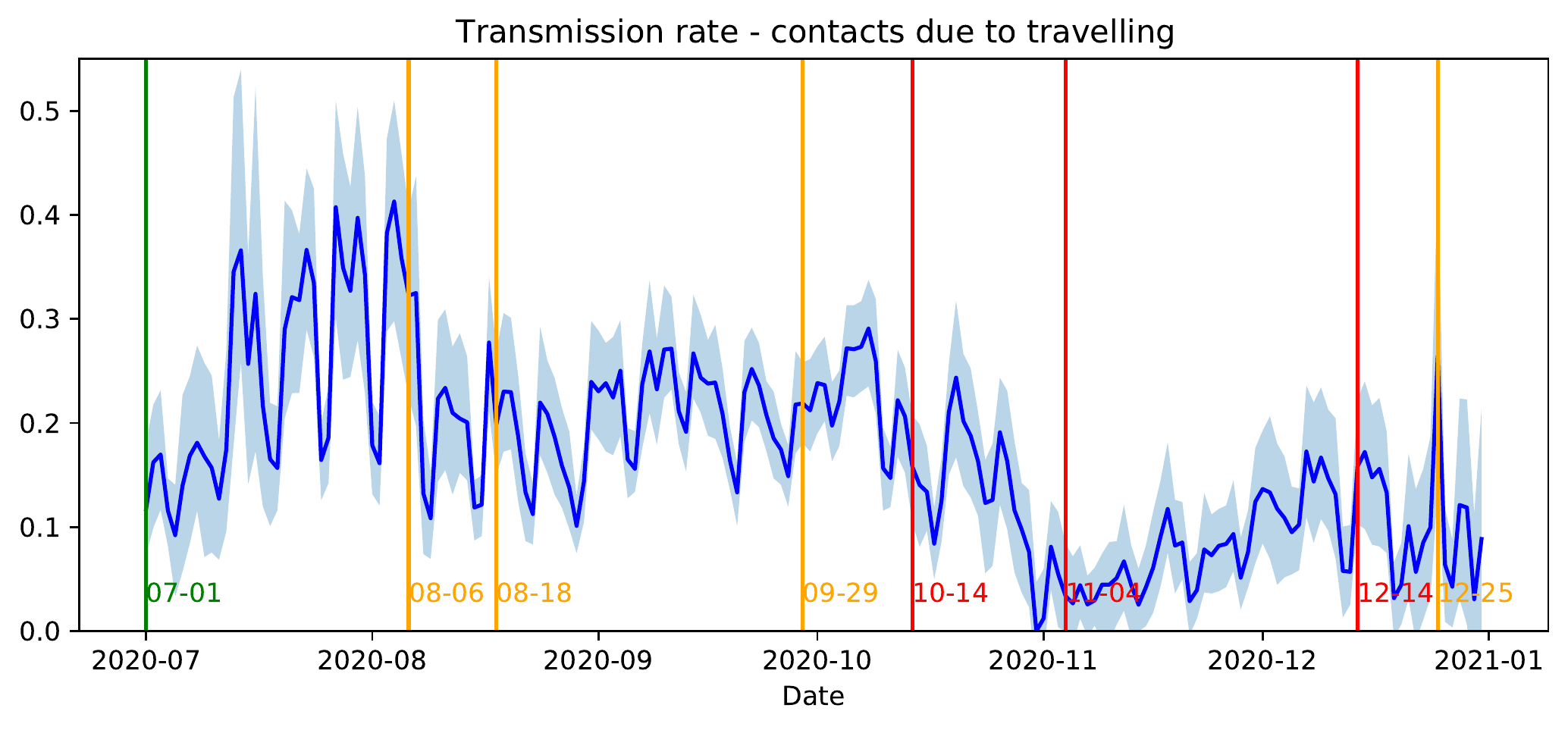}
\caption{Estimations of the mobility-induced transmission rate $\beta_{\text{mob}}$ with 95\% confidence interval.}
\label{fig_mob}
\end{figure}

Figures~\ref{fig_loc} and~\ref{fig_mob} show the estimates of the transmission rates $\beta_{\text{loc}}$ and $\beta_{\text{mob}}$. The downward peaks in $\beta_{\text{mob}}$ correspond to weekend days where the volume of commuting traffic is significantly lower than on weekdays. Generally, the behavior of the transmission rates matches in the sense that they increase and decrease in parallel. A notable exception to this can be observed in the period August 18 - September 29, where $\beta_{\text{loc}}$ increases steeply and $\beta_{\text{mob}}$ remains relatively stable. This behavior is in line with schools starting again in September after the summer break, which led to an increase in local contacts among school children and their teachers and parents.


\begin{figure}[ht]
    \centering
    \includegraphics[scale = .7]{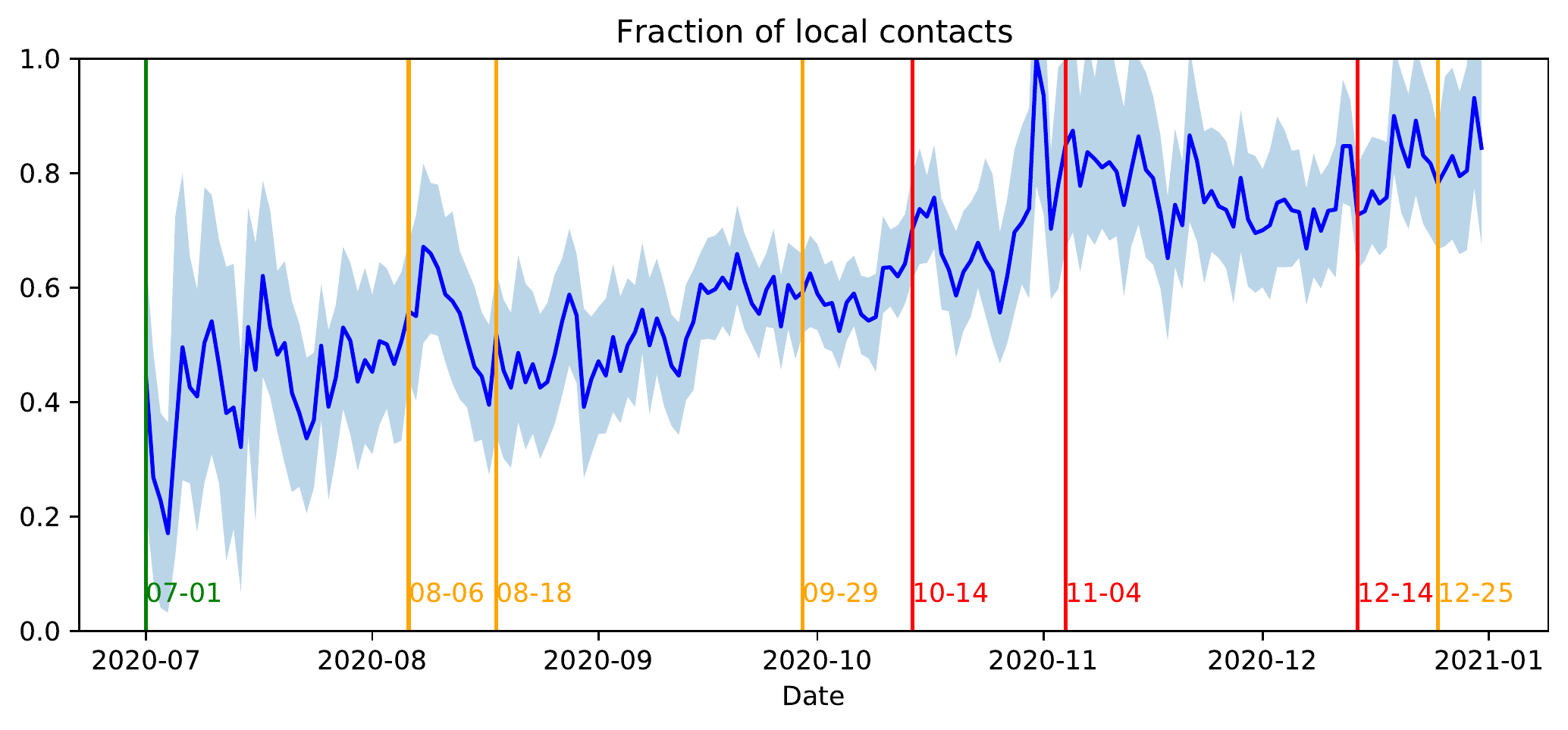}
    \caption{Estimate of the fraction $p$ of local contacts with 95\% confidence interval.}
    \label{fig:fracpos}
\end{figure}

Figure \ref{fig:fracpos} shows the estimated values for the fraction $p$ of local contacts. This fraction appears to increase steadily throughout the considered time period, which is in line with the increasing degree of mobility-restricting measures. In particular, a sudden increase occurs shortly after the partial lockdown restrictions of November 4. This suggests that our estimates capture the logical trend with respect to these measures.

\begin{figure}[ht!]
\centering
\includegraphics[scale = .7]{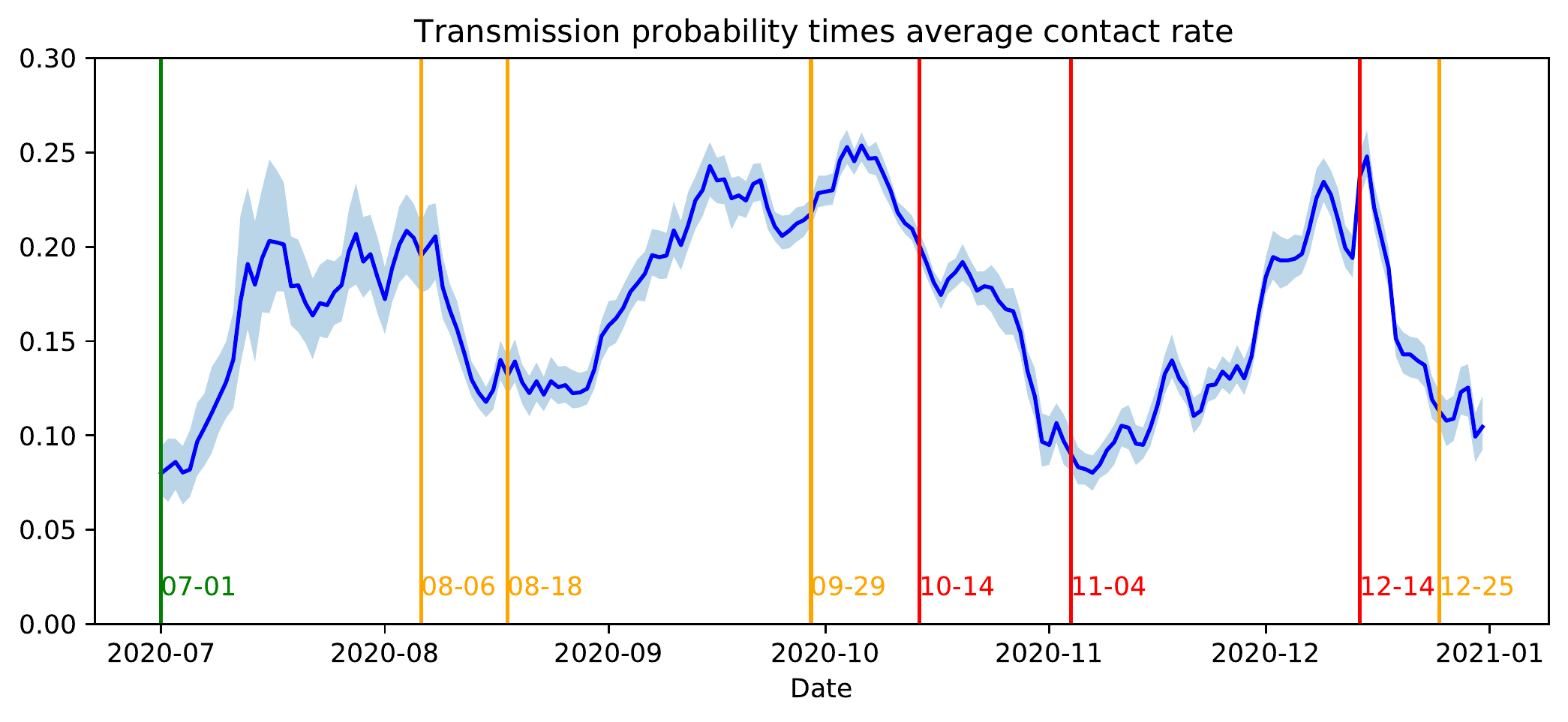}
\caption{Estimate of the term $\varepsilon c$ (transmission probability times average contact rate) with 95\% confidence interval.}
\label{fig_newlyinfected}
\end{figure}

Figure~\ref{fig_newlyinfected} shows the estimate of the term $\varepsilon c$, i.e., the transmission probability times the average contact rate. Since these two parameters appear together in the expressions for both transmission rates (see Equations (\ref{eq_beta_loc}) and (\ref{eq_beta_mob})), they cannot be estimated separately when none of them is given on forehand. However, assuming that the transmission probability remains constant over time, the estimated term $\varepsilon c$ does provide information on changes in the number of contacts. For instance, it follows from Figure~\ref{fig_newlyinfected} that the number of contacts has increased following the relaxation of measures after July 1, following the start of the new school year in September, and following the partial lockdown measures after November 4. Moreover, after the restrictions on August 6 and the (partial) lockdown measures of October 14 and December 14, the number of contacts has decreased. One explanation for this is that in this period, the government issued a strong advice to work from home whenever possible. This led to a reduction in commuting travel, which is the type of inter-regional mobility that we consider in this study.

\subsection{Validation results}

\begin{figure}[ht!]
    \centering
    \includegraphics[scale = .7]{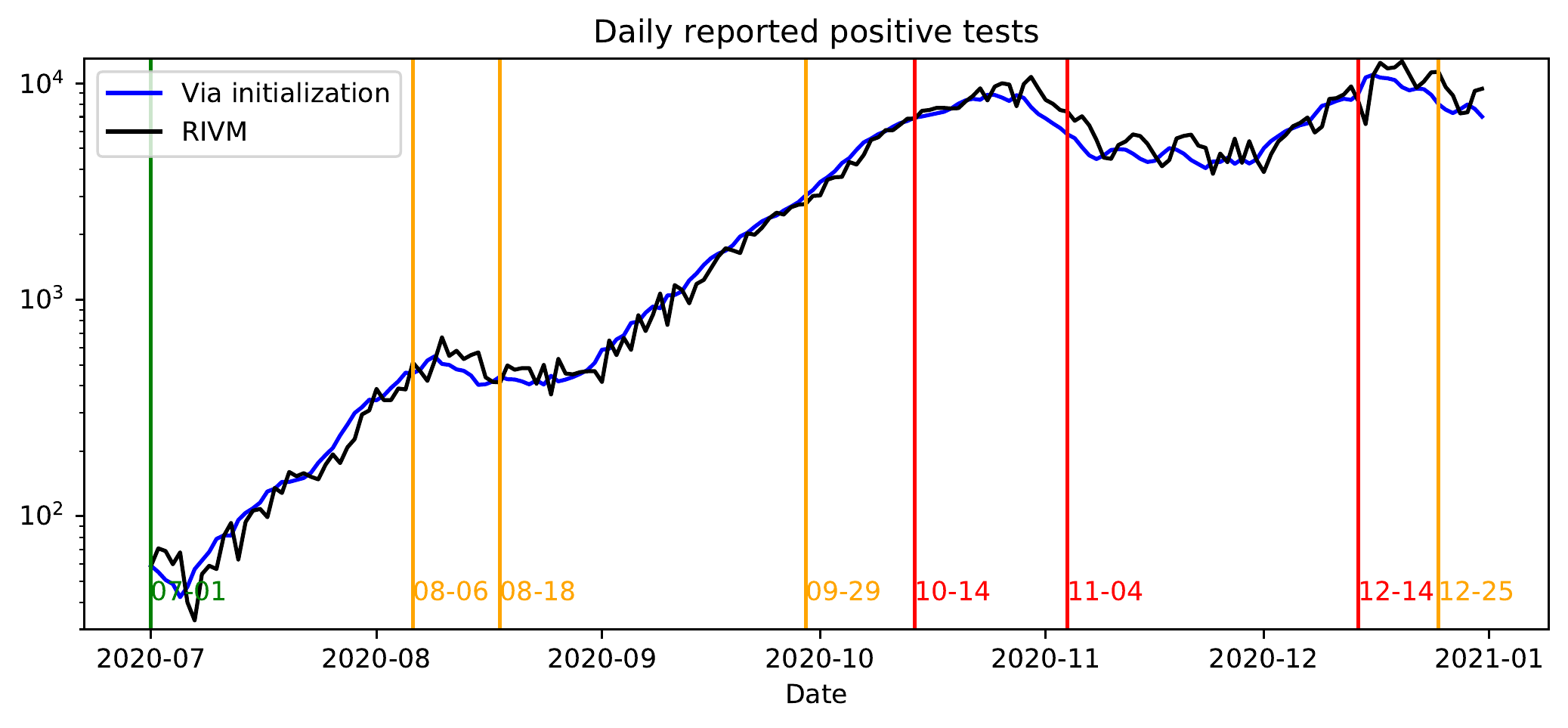}
    \caption{Validation of the initialization method.}
    \label{fig_final_newtests_init}
\end{figure}

\begin{figure}[ht!]
    \centering
    \includegraphics[scale = .7]{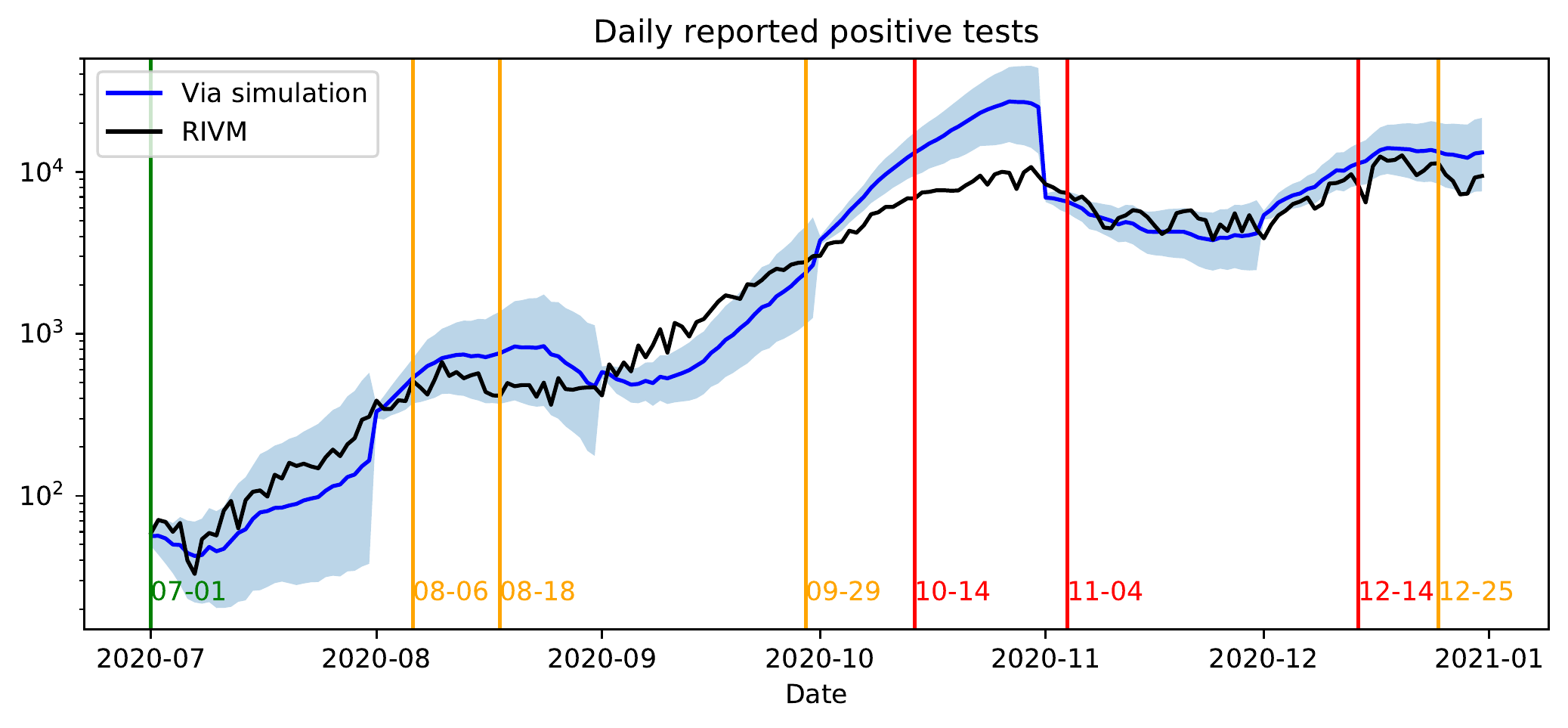}
    \caption{Validation of the simulation method with 95\% confidence intervals.}
    \label{fig_final_newtests_sim}
\end{figure}

Figures~\ref{fig_final_newtests_init} and~\ref{fig_final_newtests_sim} show the estimates of the daily reported positive tests computed via the initialization and simulation method, respectively (see Section~\ref{sec_val_nat}), and those as reported by RIVM. Figure~\ref{fig_final_newtests_init} indicates that the estimation via the initialization method closely follows the positive tests as reported by RIVM. This suggests that our initialization approach is able to estimate the state of the epidemic well on the national level. Figure~\ref{fig_final_newtests_sim} indicates that the simulation method is generally able to follow the general trend of the reported new tests. This suggests that our model and our methods for parameter initialization and estimation are successful in capturing the development of the number of daily reported positive tests on the national level. Moreover, our results suggest that the used data sources on commuting mobility and SARS-CoV-2 seroprevalence can be integrated successfully into compartmental models.

To determine whether our estimates of the effective reproduction number have a time delay compared to those of RIVM, we computed the correlation between these two estimates for different numbers of delay of days. These results show that there is a strong correlation between these estimates when the delay is between 5-14 days and is the strongest for a delay of 11 days ($0.865$). This is in line with the reporting delay and difference in time window and data sources between our method and the methods used by RIVM.

\begin{figure}[ht!]
    \centering
    \includegraphics[scale = .7]{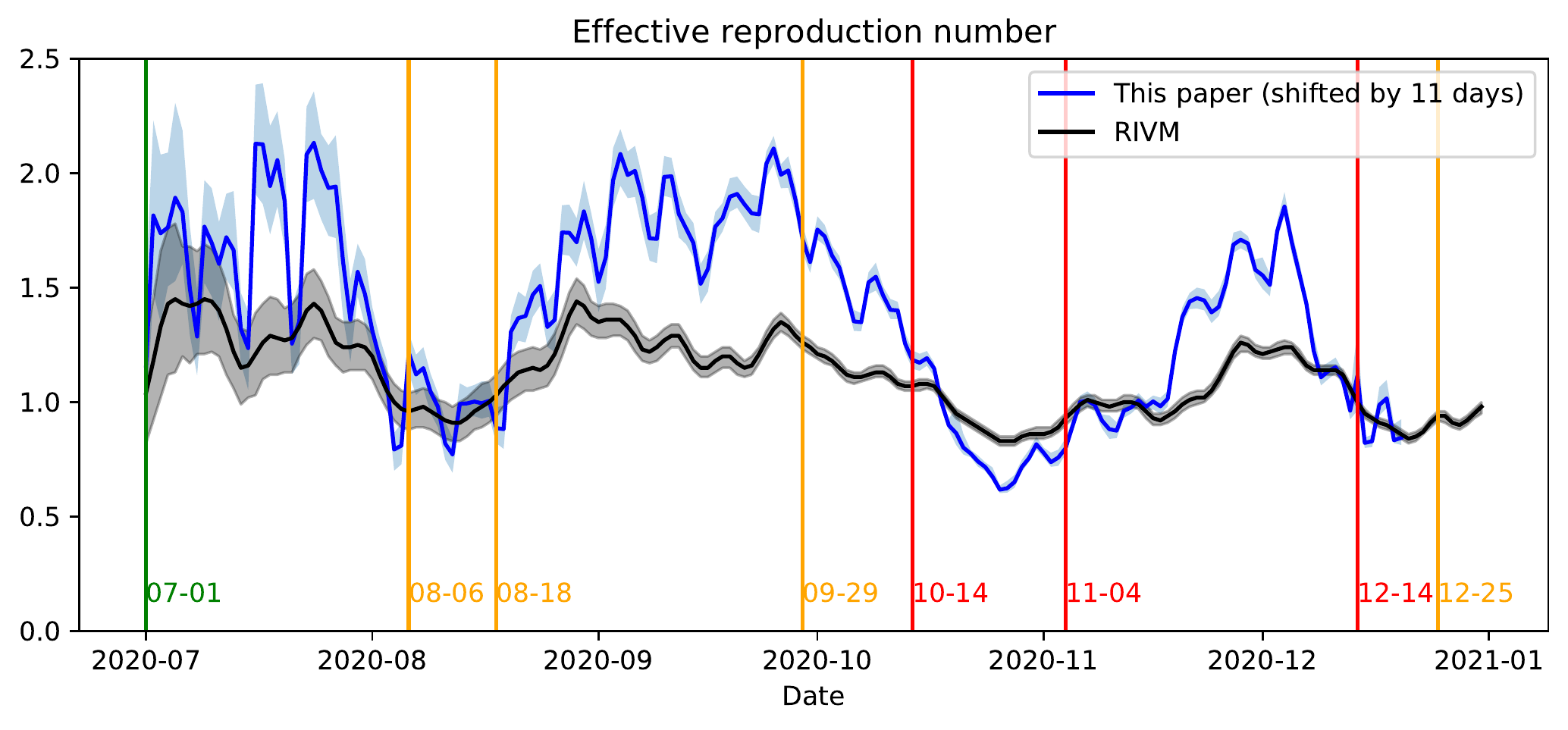}
    \caption{The effective reproduction number as estimated in this paper (shifted by 11 days) and by RIVM with 95\% confidence intervals.}
    \label{fig_final_R_shift}
\end{figure}

Figure~\ref{fig_final_R_shift} shows the estimates for the effective reproduction number, shifted by 11 days, and the estimates made by RIVM. Overall, our estimates appear to be somewhat more extreme than those of RIVM. However, both estimates share the same general pattern with regard to periods of monotonic increases and decreases and agree on whether the epidemic is expanding or dying out, i.e., whether the effective reproduction number is larger or smaller than 1. Moreover, the peaks and valleys in the estimates match very well, which confirms the relative consistency of the estimates and the presence of a temporal shift. Overall, the validation shows that our estimation method agrees qualitatively with the method used by RIVM.

\subsection{Forecasting accuracy}

\begin{figure}[ht!]
    \centering
    \includegraphics[scale = .7]{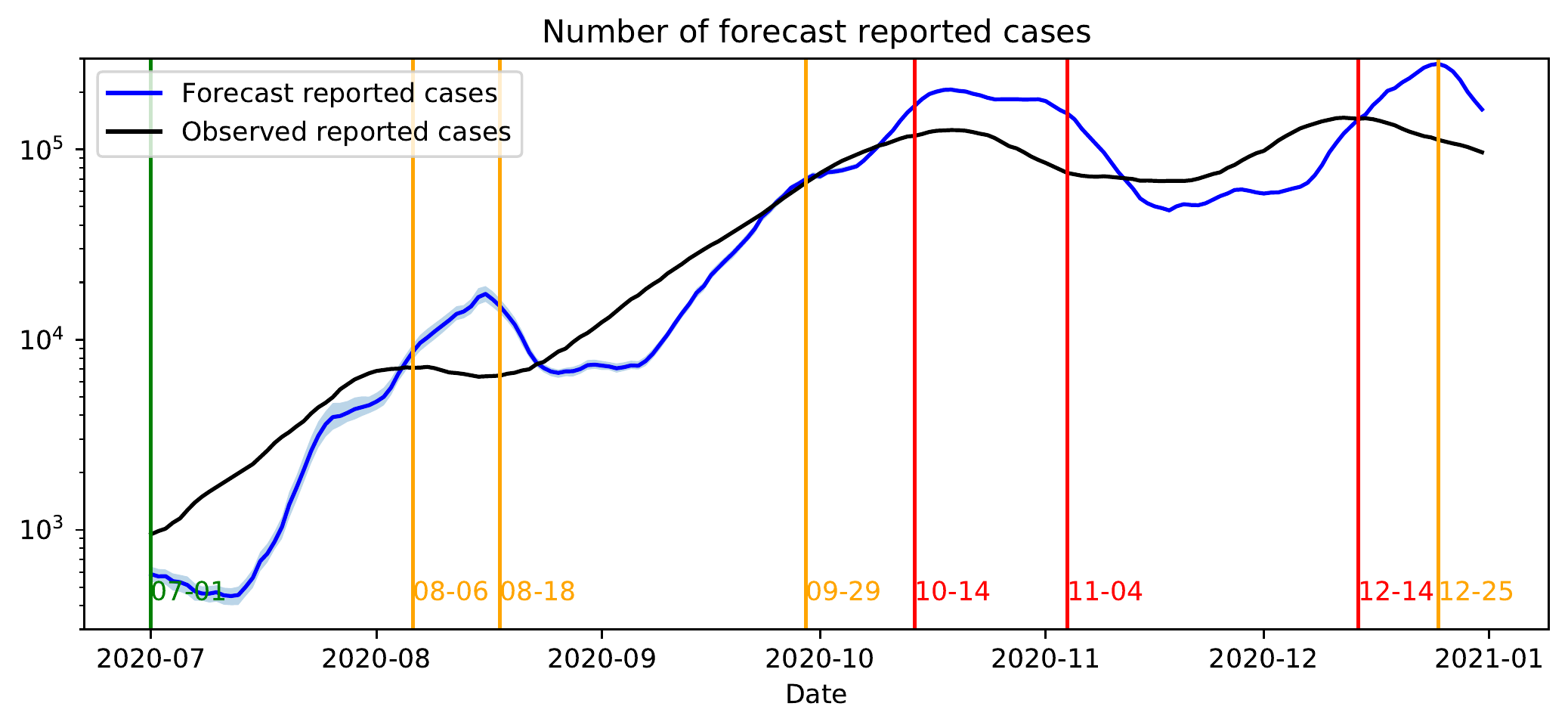}
    \caption{Forecast and actual number of reported cases in the following 14 days with 95\% confidence interval.}
    \label{fig_final_tot_14}
\end{figure}

In this subsection, we assess the forecasting capabilities of the compartmental model when parameters are estimated as described in Section~\ref{sec_setup_forecast}. First, we focus on the accuracy of our forecasts of the \emph{absolute} number of daily reported cases. Figure~\ref{fig_final_tot_14} shows for each date within the period July 1 - December 31, 2020, our forecast and the actual number of reported cases nation-wide in the following 14 days. This figure shows that, generally, the forecasts follow the actual number of reported cases. Analogous to the estimation of the effective reproduction number, it appears that the forecasts have a certain delay. One reason for this is that our model is necessarily initialized using parameters from 1 to 2 weeks ago as input. When the transmission rates are relatively stable for a longer period of time, this effect is minimal. In this case, the forecast error is relatively small (e.g., around August 1 and in October). However, when the transmission rates are not stable, for example due to newly introduced preventive measures, the forecasts are made using data that is not representative of the current epidemiological situation. Consequently, it takes some time before the input data of the forecast procedure is again a proper representation of the future parameter values (see also the behavior of the transmission rates in Figures~\ref{fig_loc} and~\ref{fig_mob}).

\begin{figure}
    \centering
    \begin{subfigure}[t]{.49\textwidth}
    \centering
        \includegraphics[scale = .35]{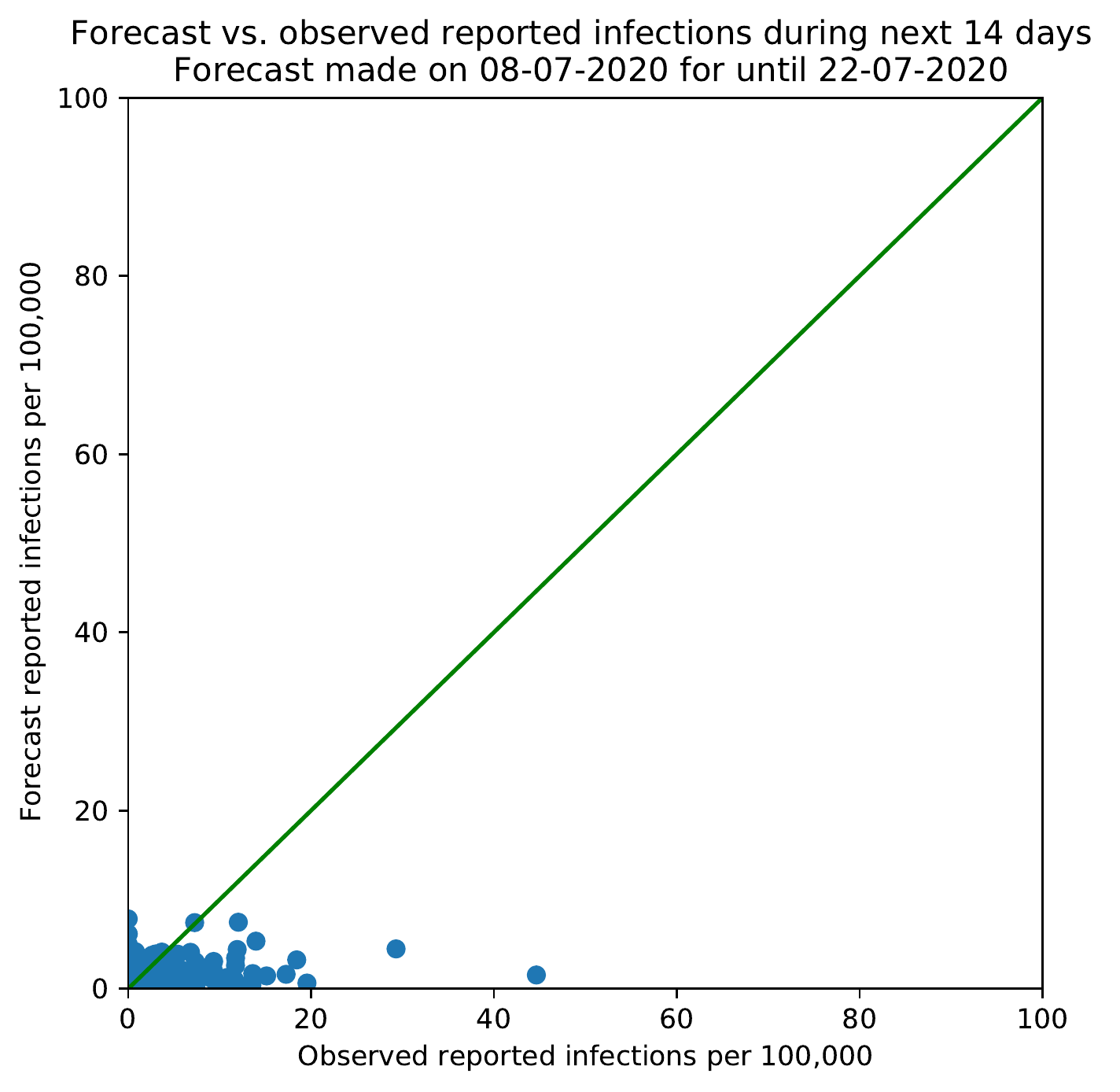}
    \end{subfigure}
        \begin{subfigure}[t]{.49\textwidth}
        \centering
        \includegraphics[scale = .35]{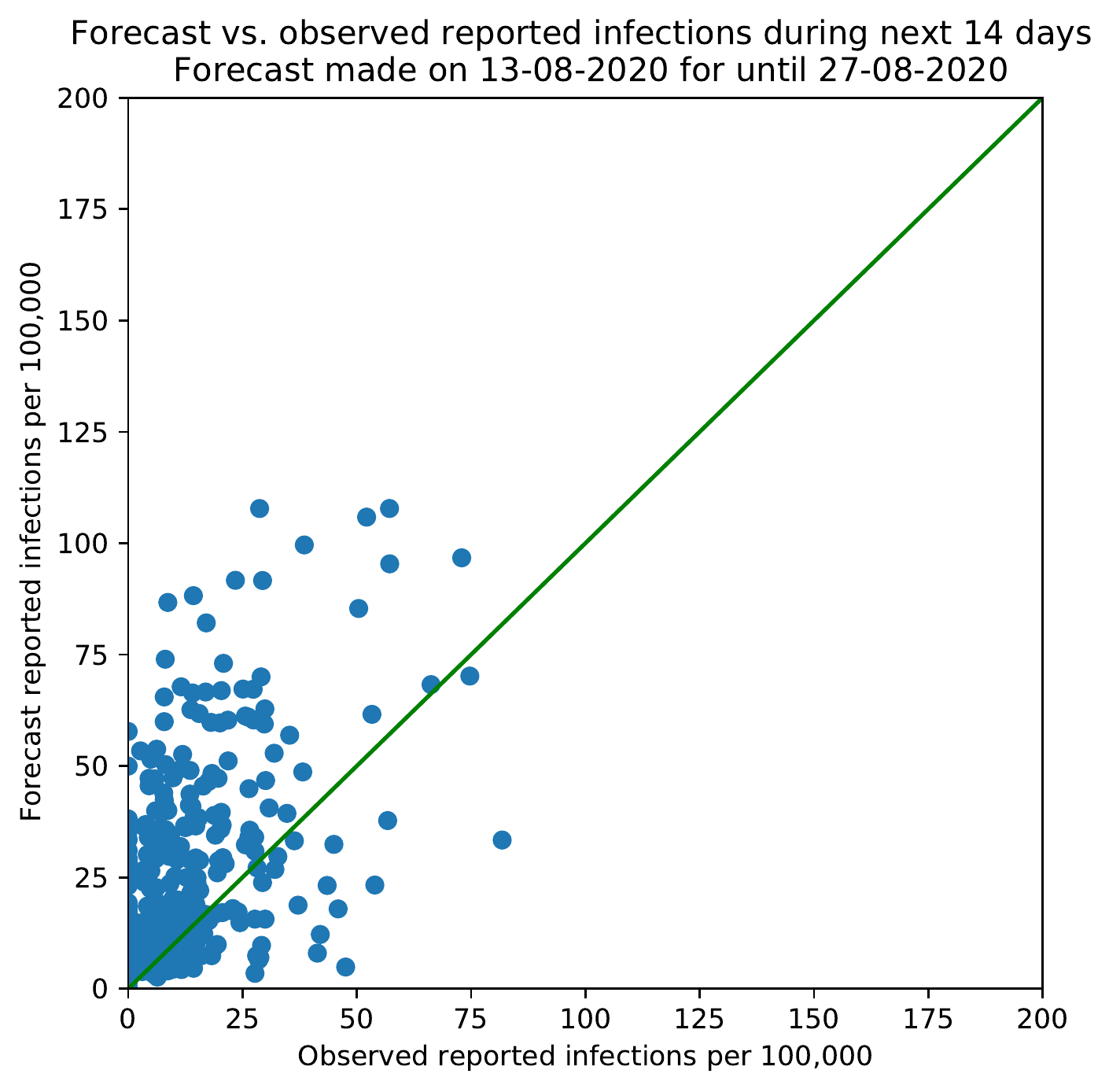}
    \end{subfigure}
     \begin{subfigure}[t]{.49\textwidth}
     \centering
        \includegraphics[scale = .35]{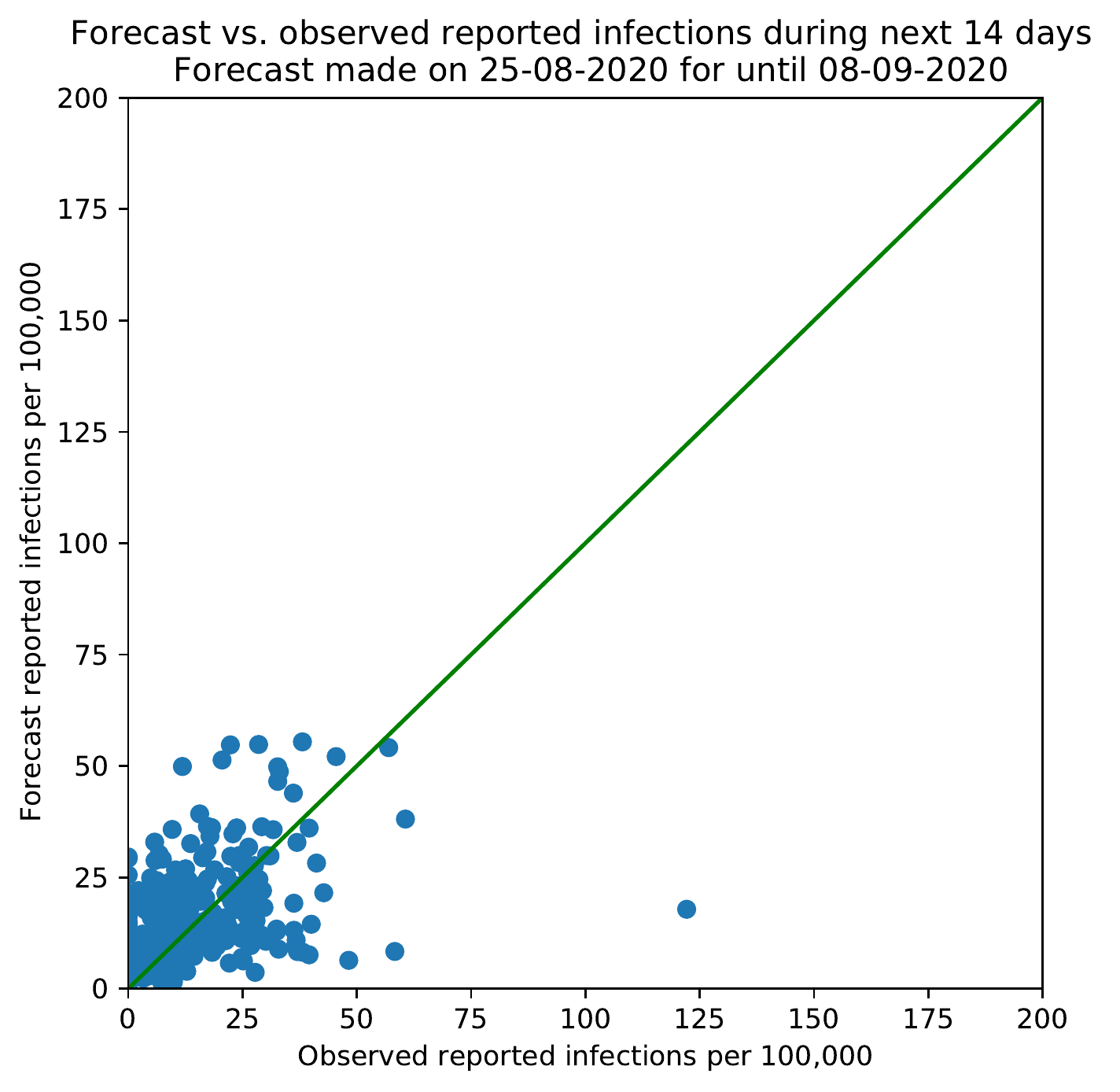}
    \end{subfigure}
        \begin{subfigure}[t]{.49\textwidth}
        \centering
        \includegraphics[scale = .35]{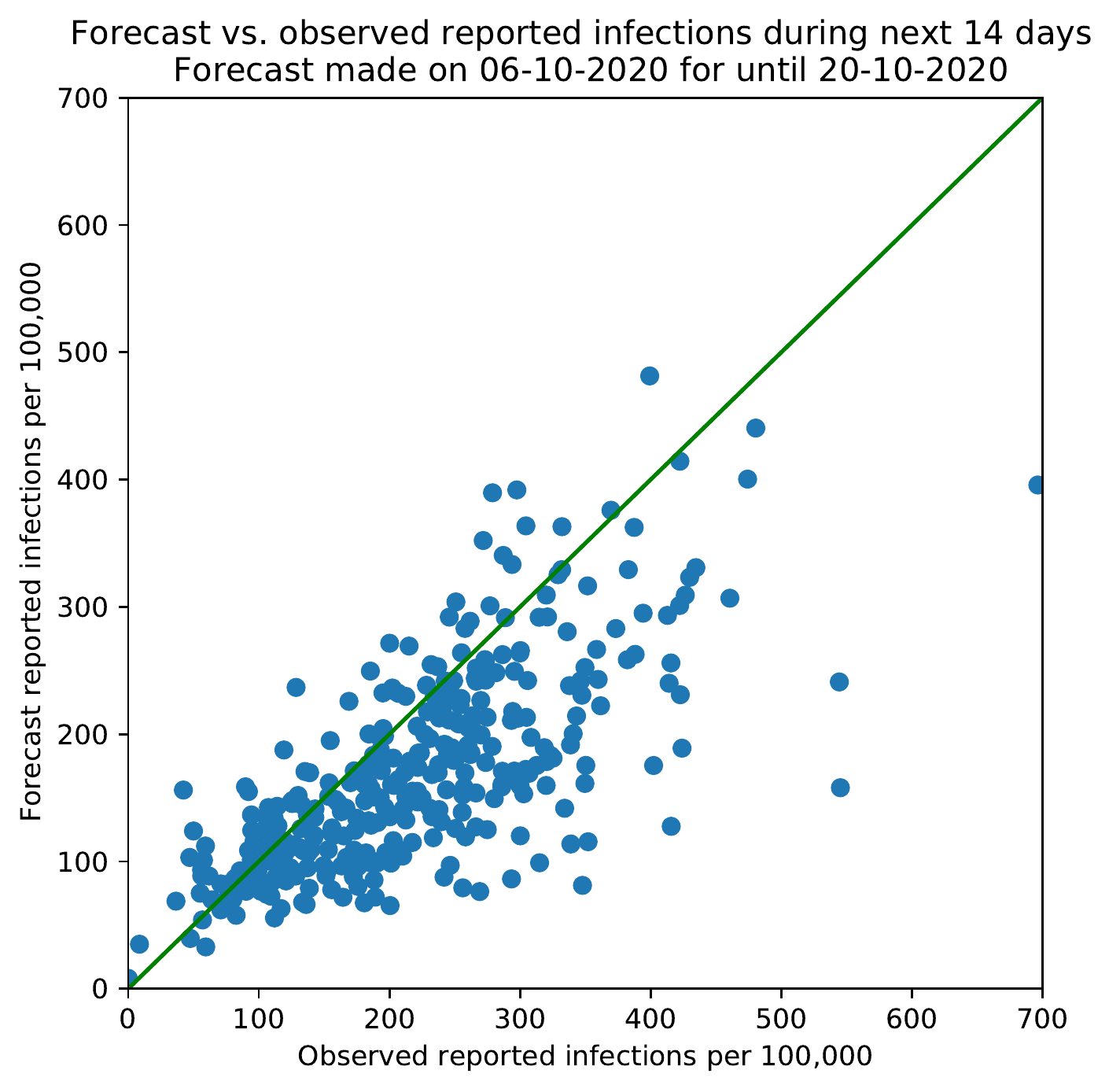}
    \end{subfigure}
        \begin{subfigure}[t]{.49\textwidth}
    \centering
        \includegraphics[scale = .35]{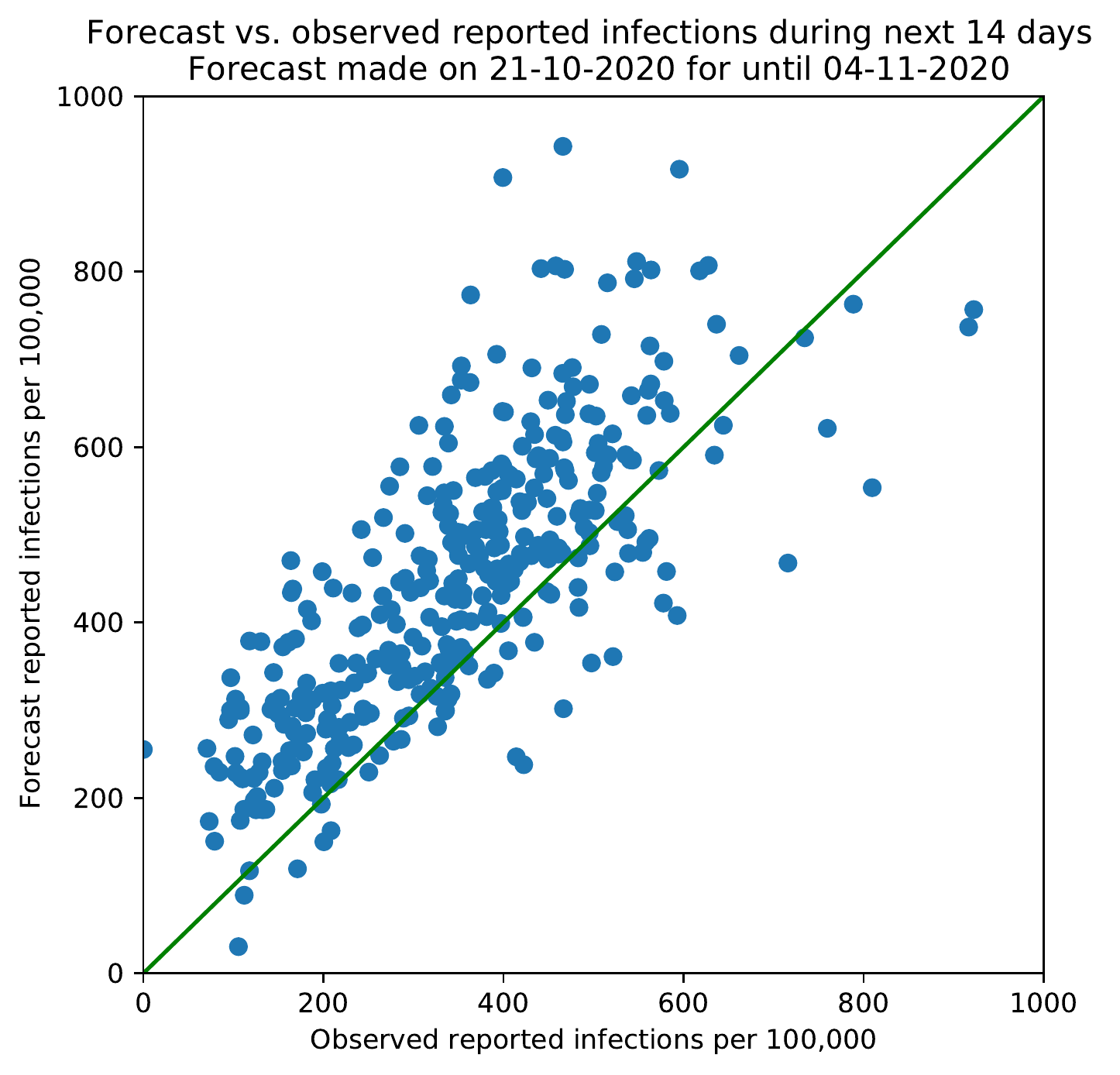}
    \end{subfigure}
        \begin{subfigure}[t]{.49\textwidth}
        \centering
        \includegraphics[scale = .35]{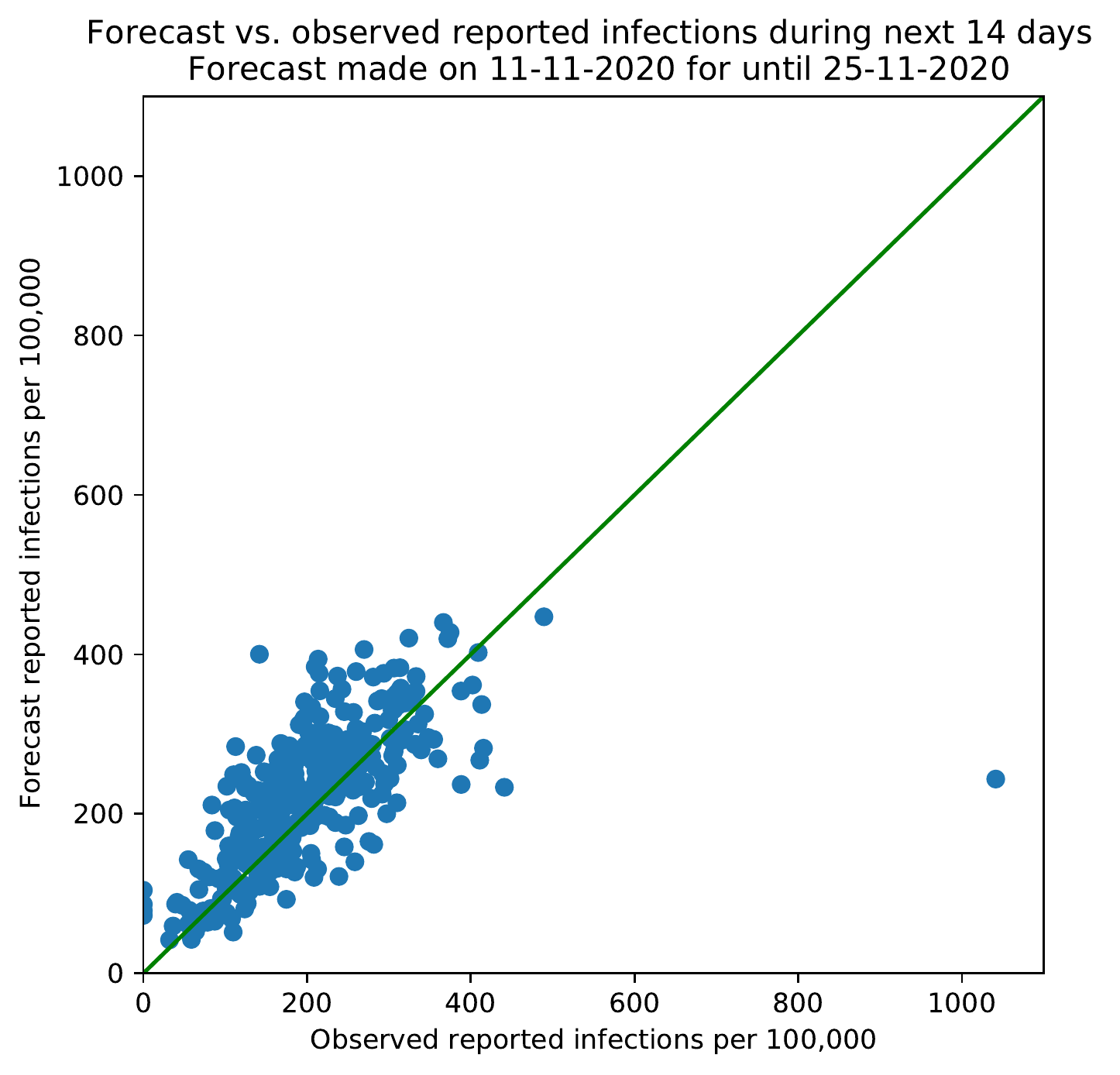}
    \end{subfigure}
     \begin{subfigure}[t]{.49\textwidth}
     \centering
        \includegraphics[scale = .35]{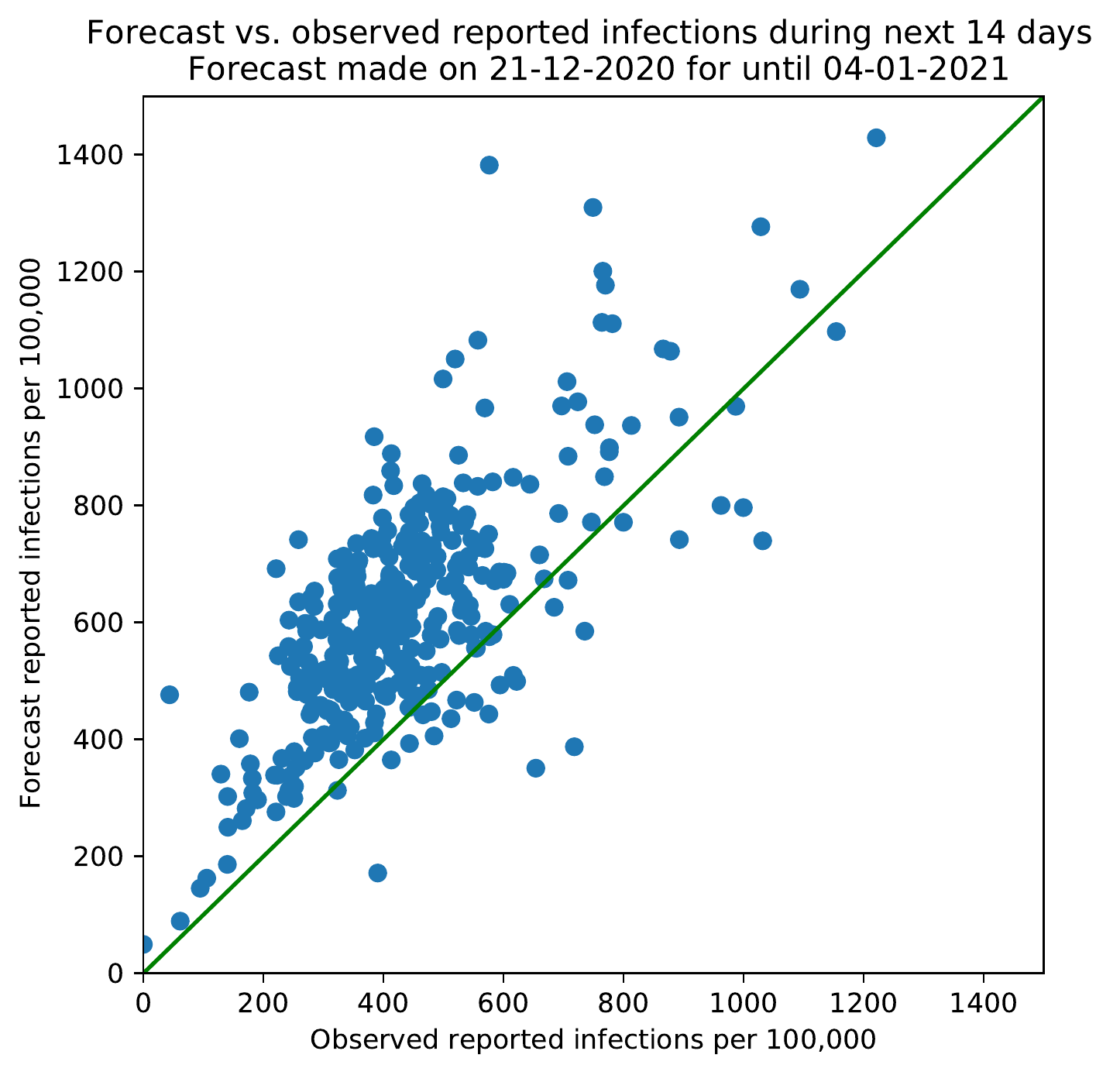}
    \end{subfigure}
        \begin{subfigure}[t]{.49\textwidth}
        \centering
        \includegraphics[scale = .35]{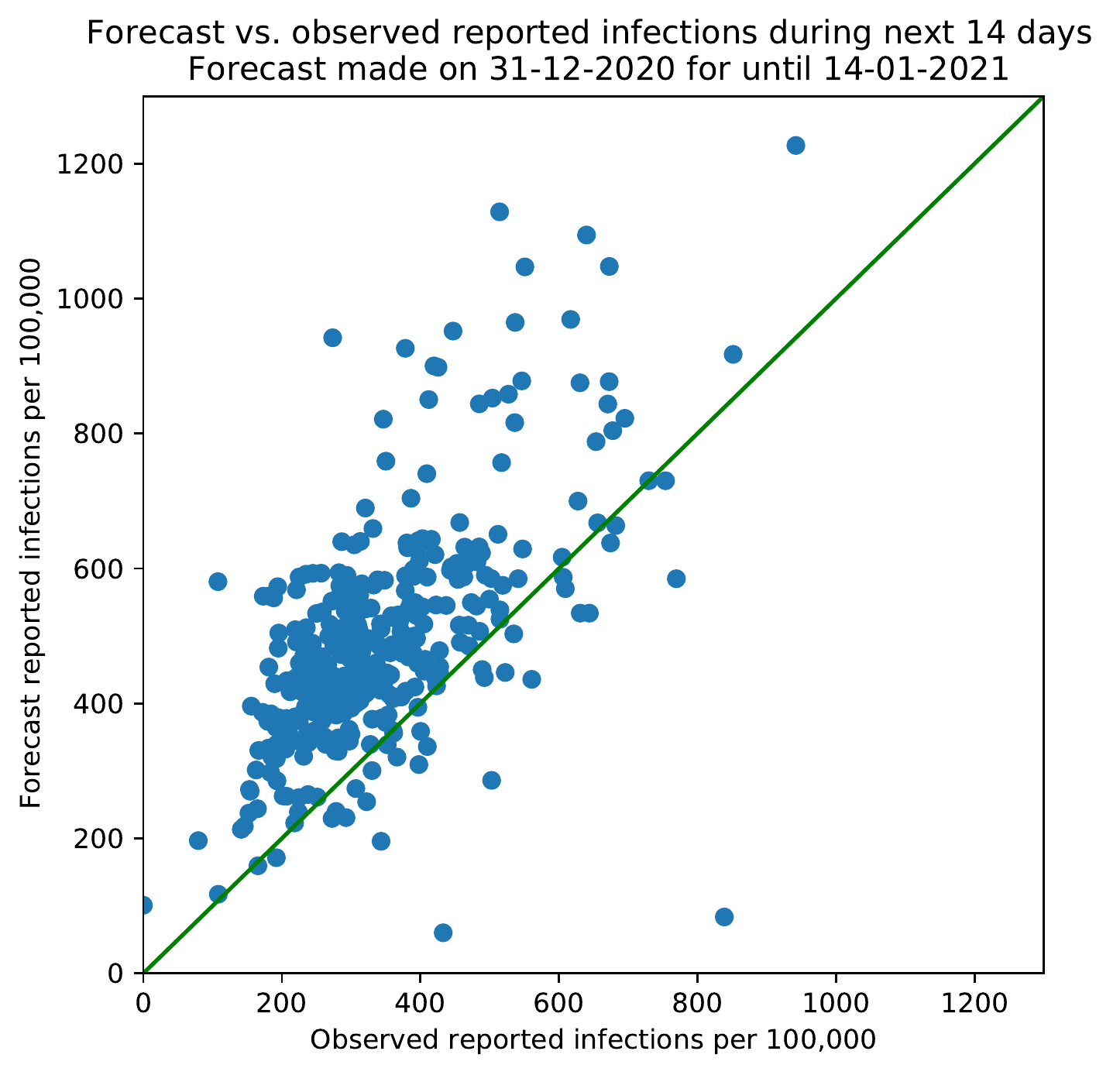}
    \end{subfigure}
    \caption{Scatter-plots of the accuracy of the 14-day ahead forecasts made on each of the eight selected dates.}
    \label{fig_scatter}
\end{figure}

\begin{figure}
    \centering
    \begin{subfigure}[t]{.49\textwidth}
    \centering
        \includegraphics[scale = .45]{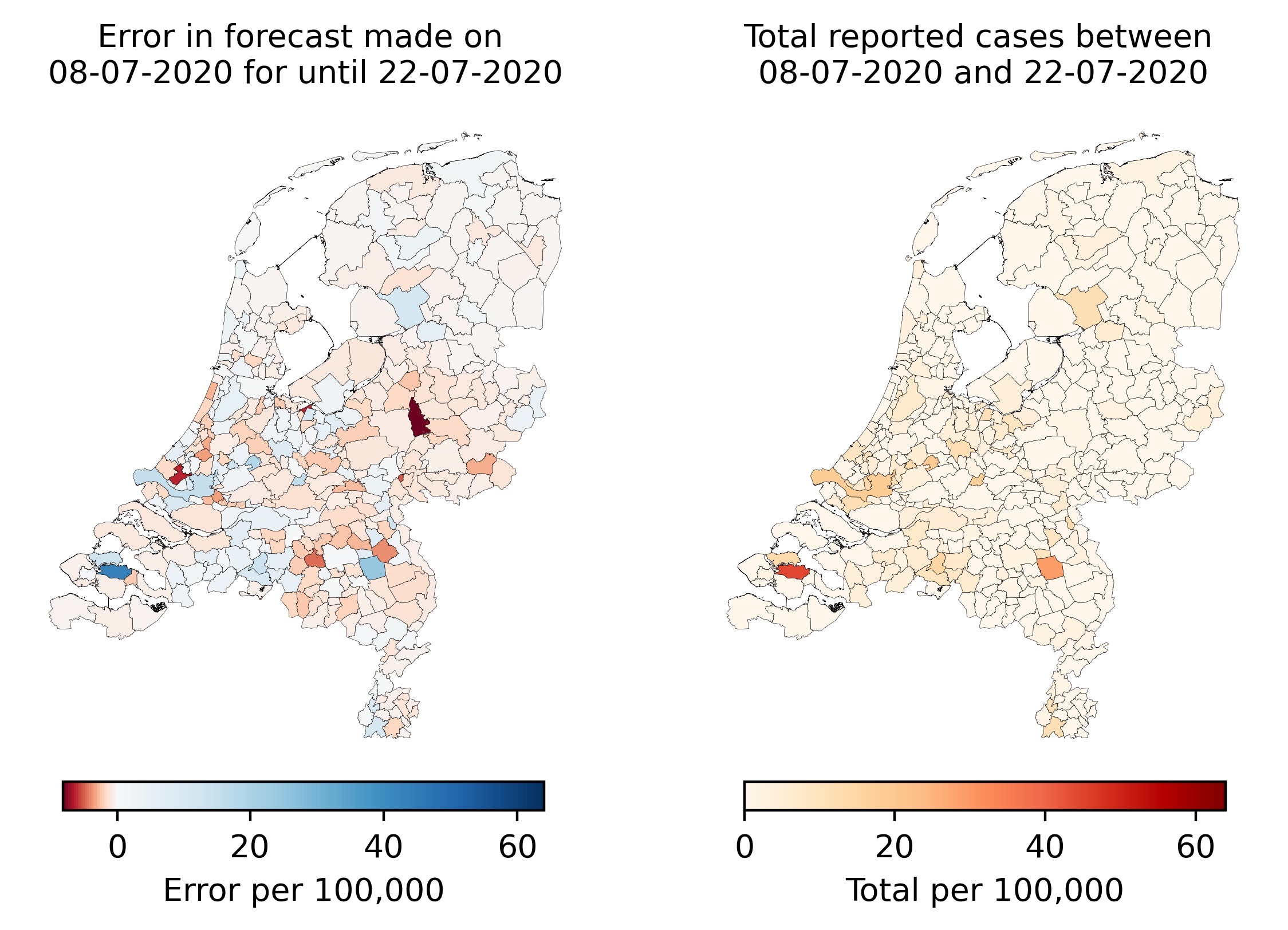}
    \end{subfigure}
        \begin{subfigure}[t]{.49\textwidth}
        \centering
        \includegraphics[scale = .45]{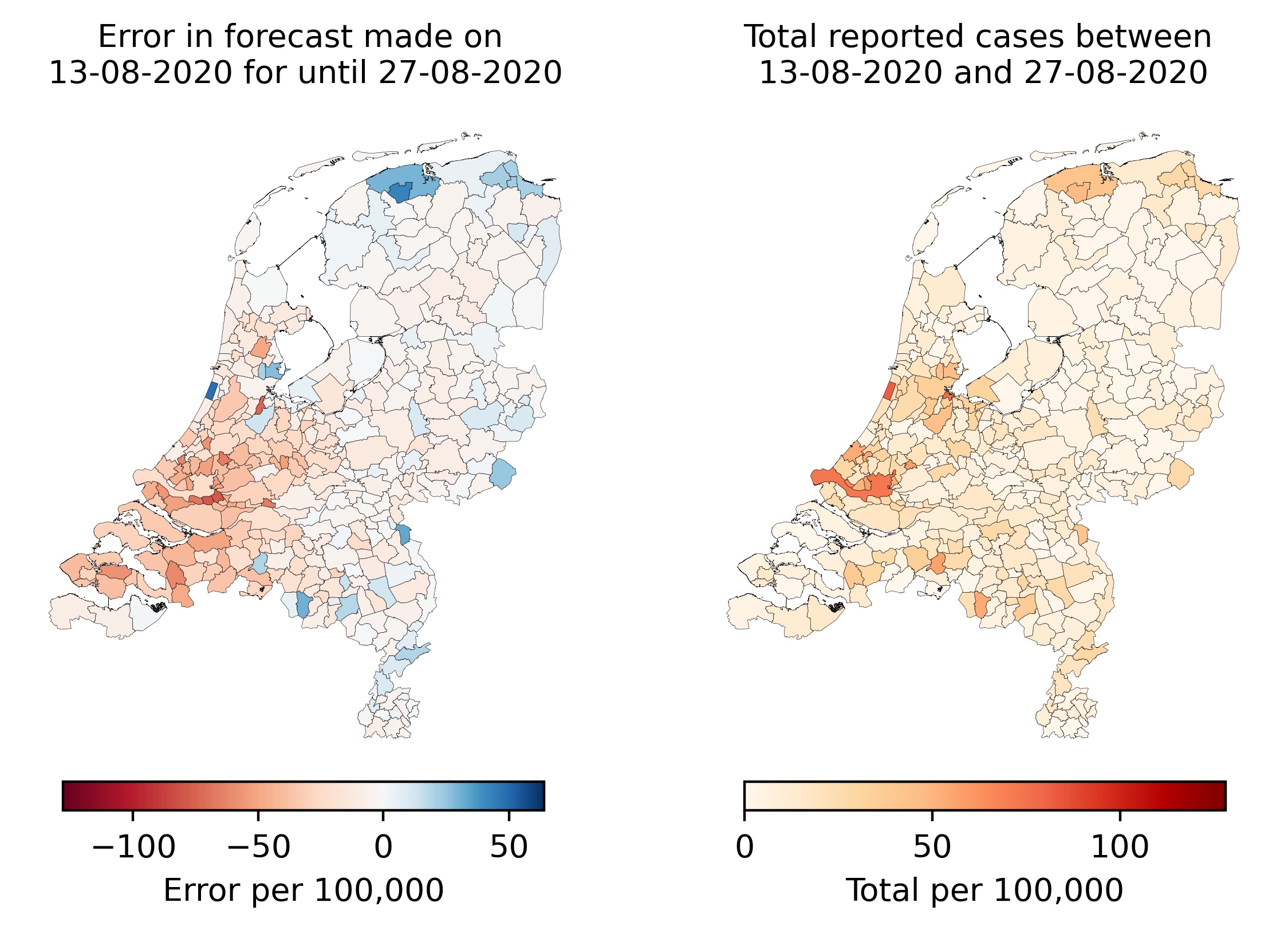}
    \end{subfigure}
     \begin{subfigure}[t]{.49\textwidth}
     \centering
        \includegraphics[scale = .45]{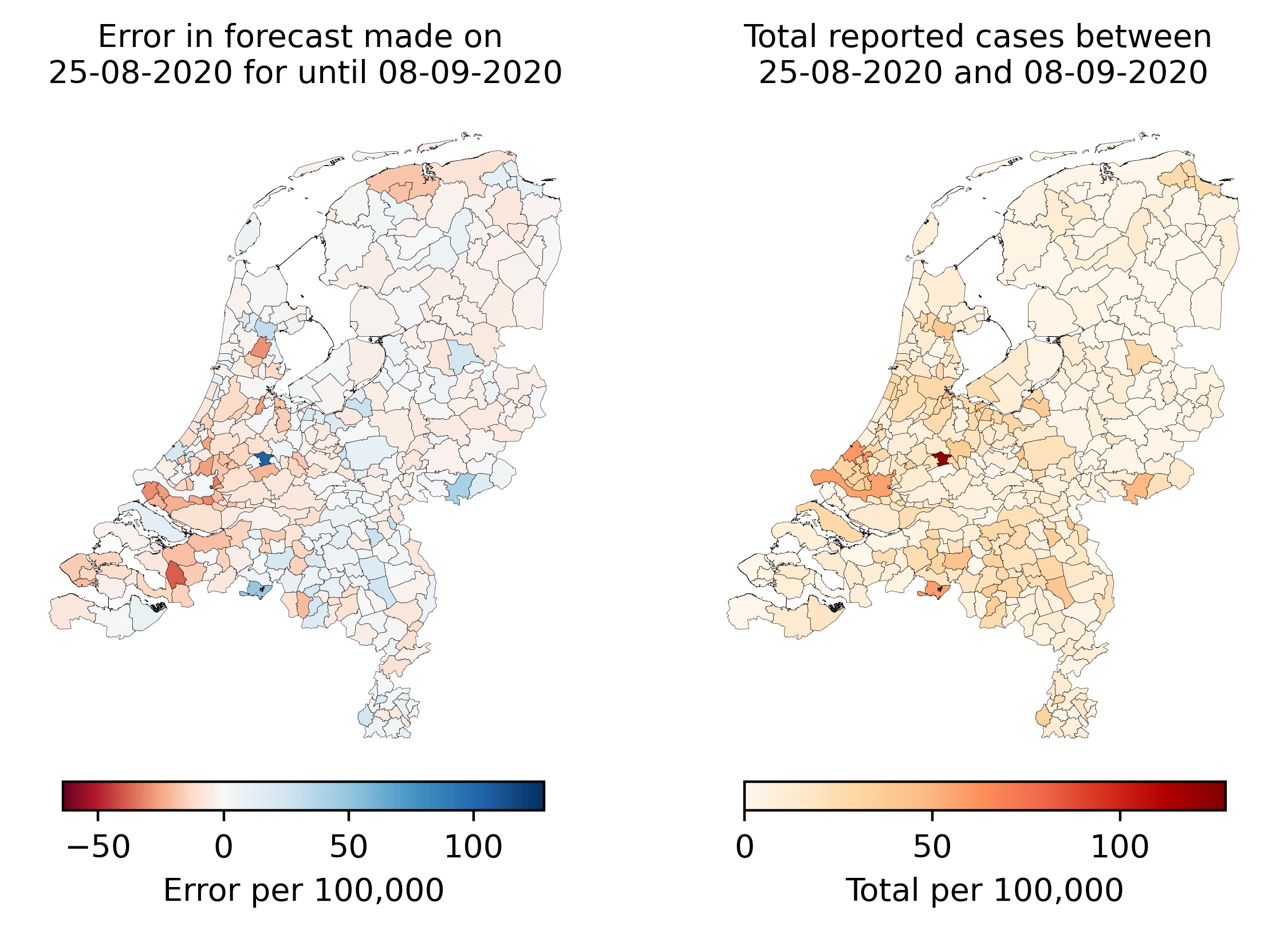}
    \end{subfigure}
        \begin{subfigure}[t]{.49\textwidth}
        \centering
        \includegraphics[scale = .45]{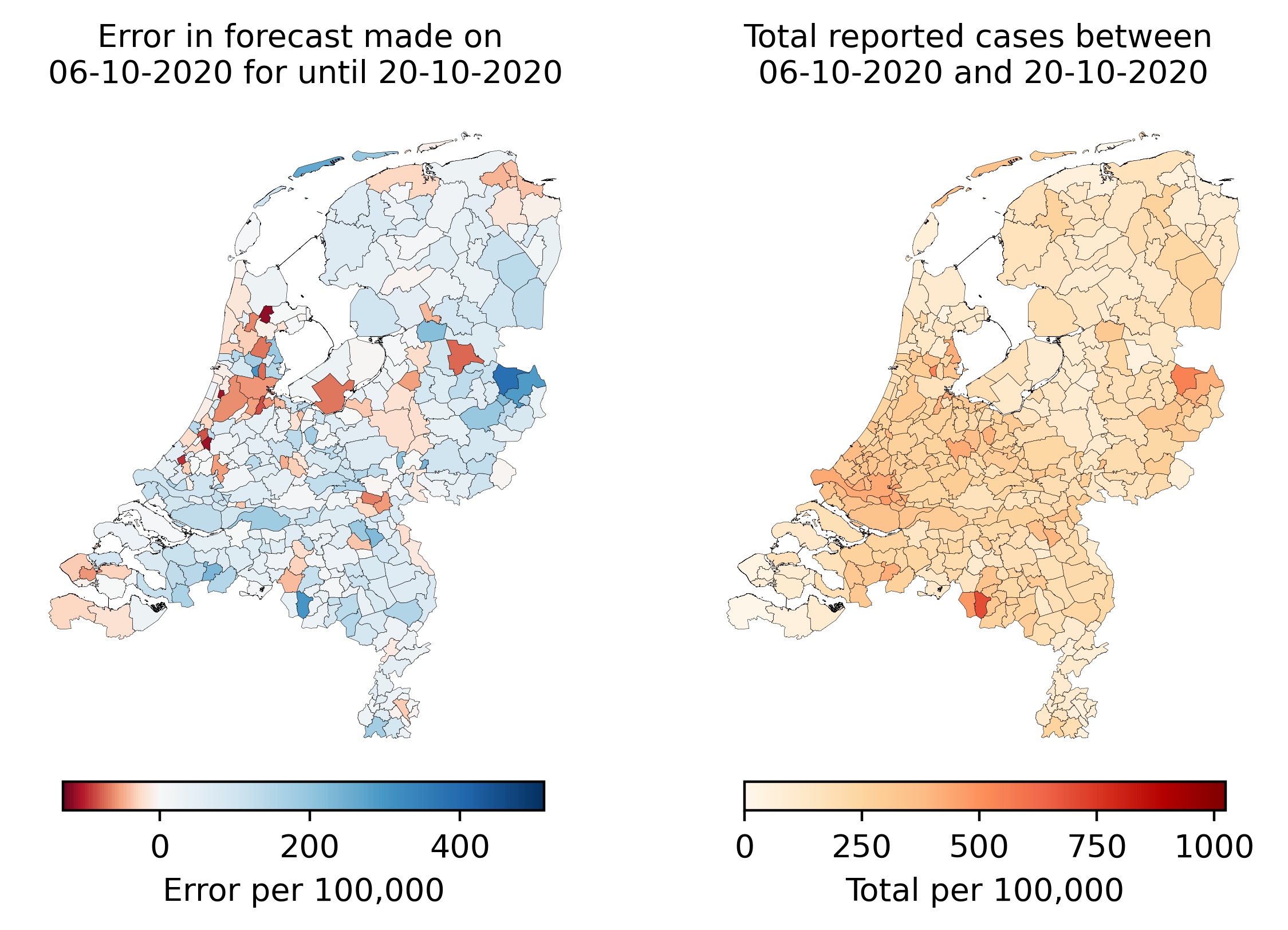}
    \end{subfigure}
        \begin{subfigure}[t]{.49\textwidth}
    \centering
        \includegraphics[scale = .45]{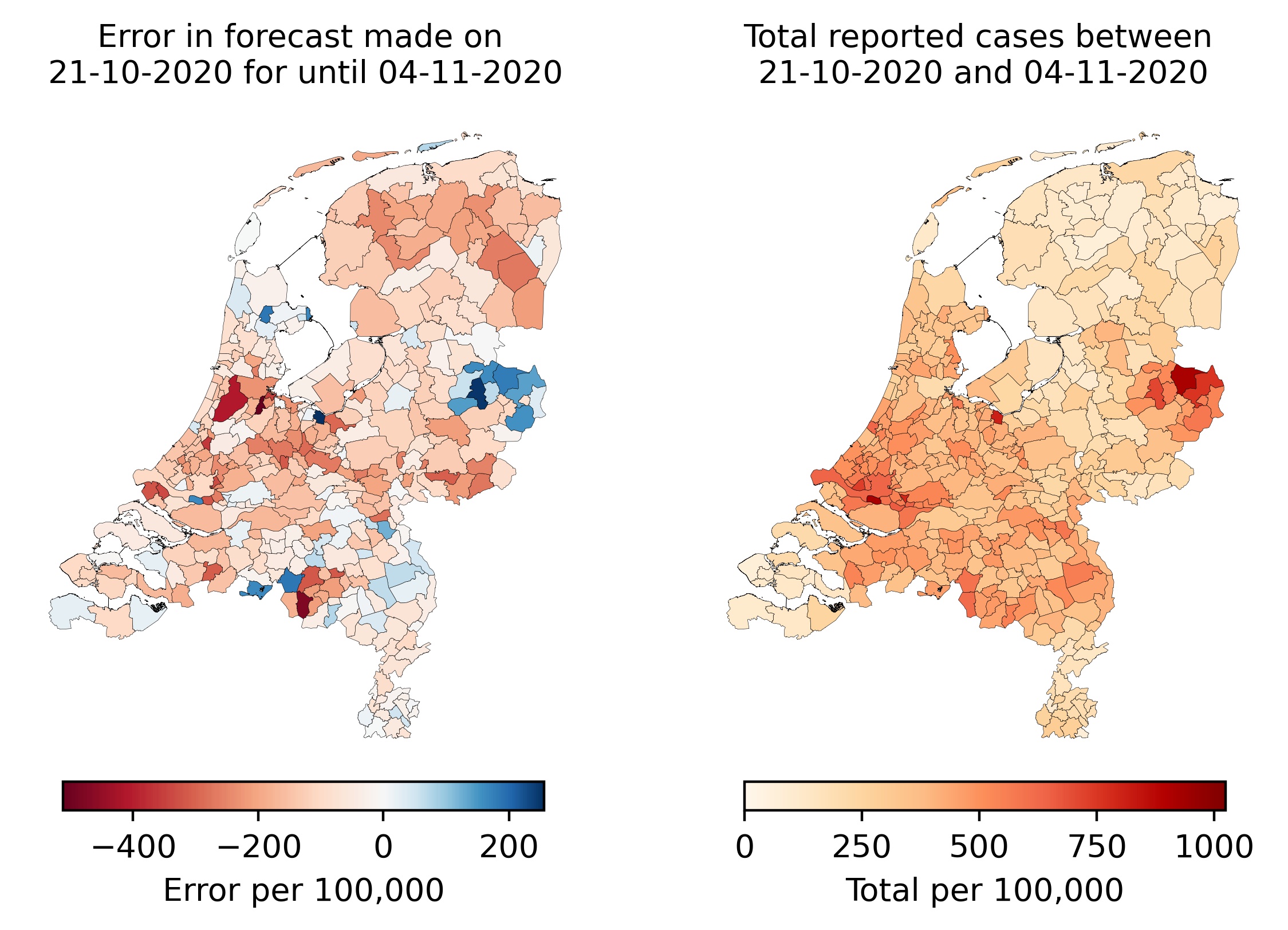}
    \end{subfigure}
        \begin{subfigure}[t]{.49\textwidth}
        \centering
        \includegraphics[scale = .45]{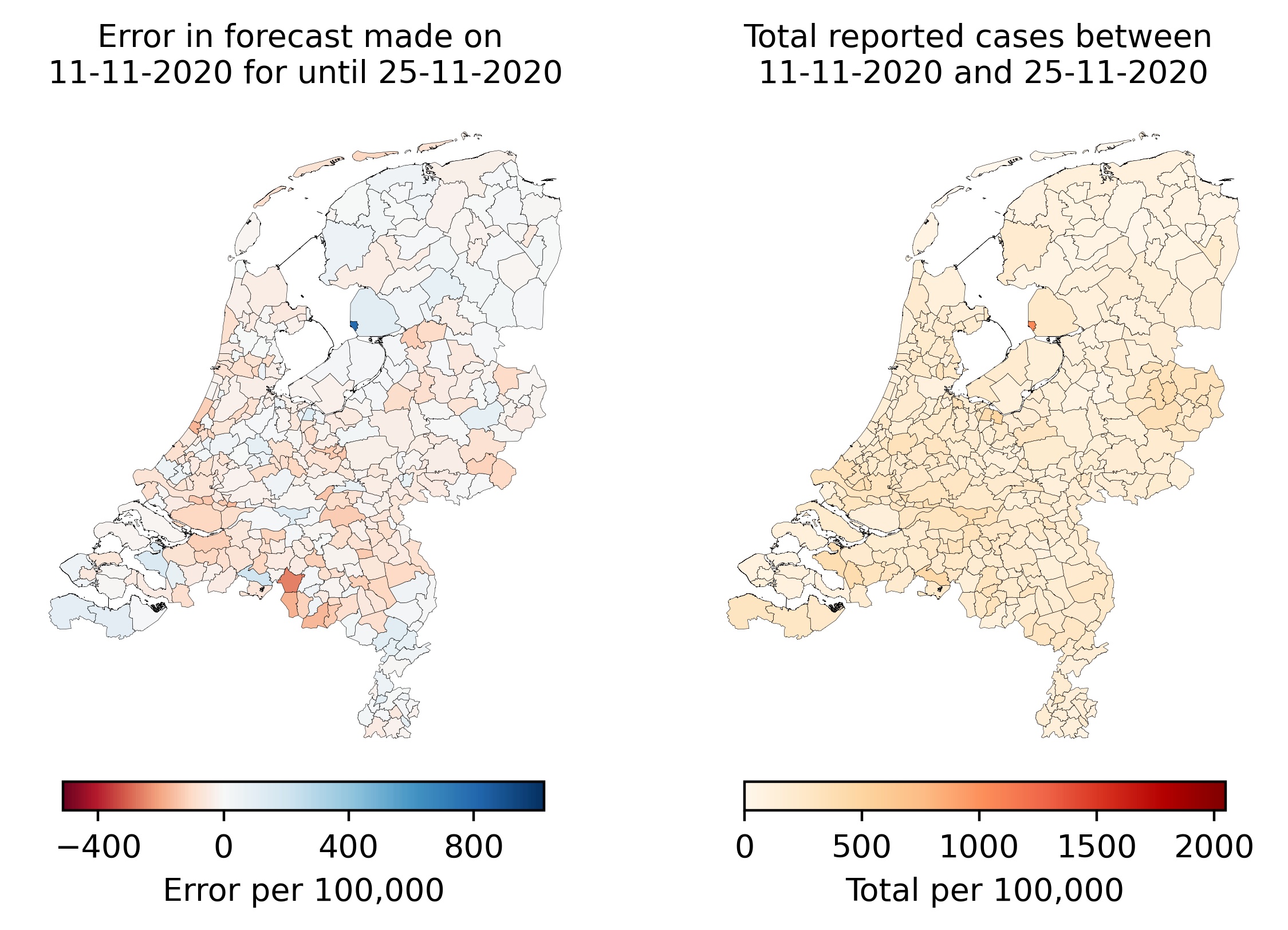}
    \end{subfigure}
     \begin{subfigure}[t]{.49\textwidth}
     \centering
        \includegraphics[scale = .45]{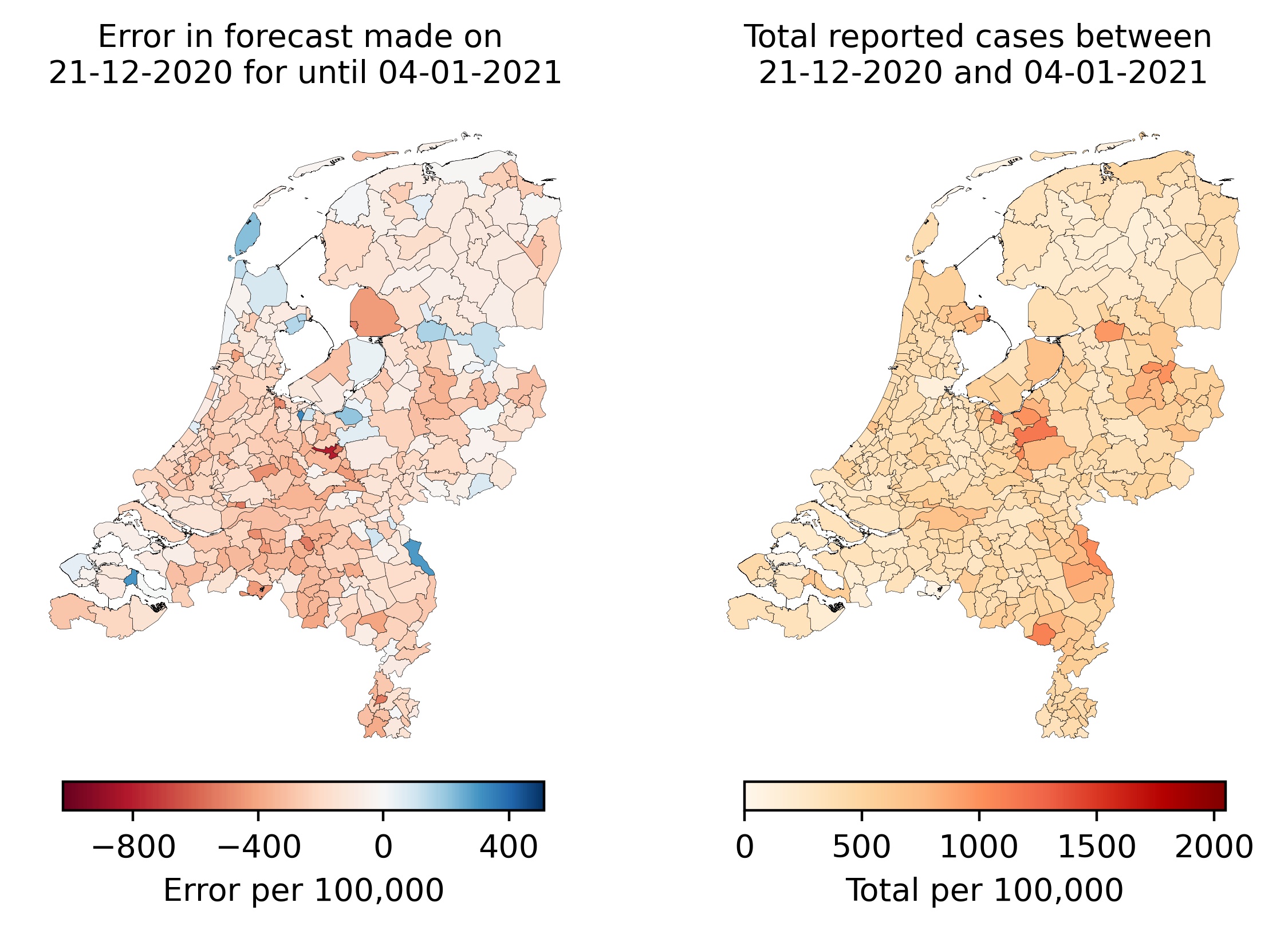}
    \end{subfigure}
        \begin{subfigure}[t]{.49\textwidth}
        \centering
        \includegraphics[scale = .45]{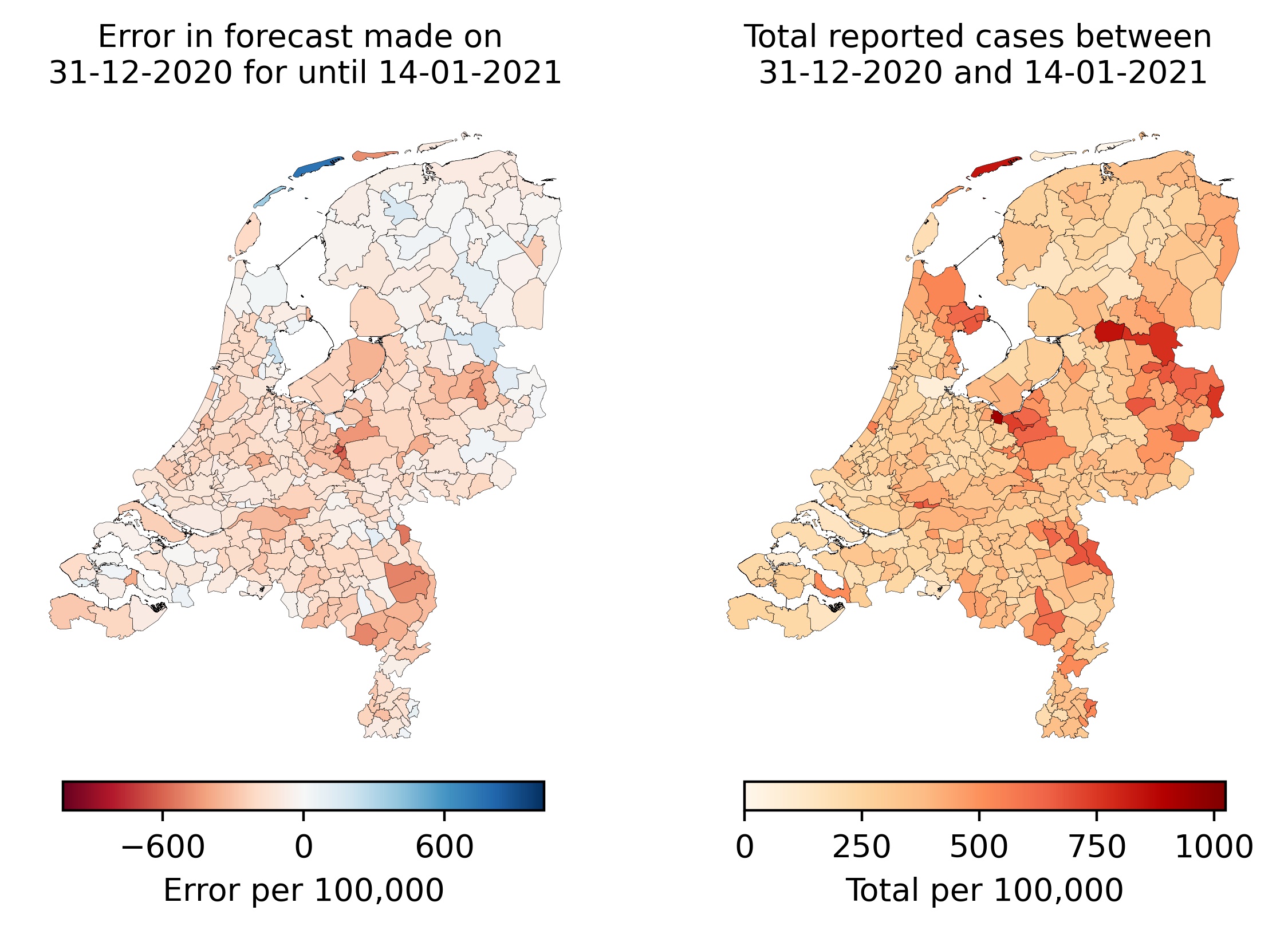}
    \end{subfigure}
    \caption{Geographical representation of the error in the 14-day ahead forecasts made on each of the eight selected dates.}
    \label{fig_maps}
\end{figure}

We now focus on the forecast accuracy at the level of municipalities. We discuss the results for eight specific dates, namely those that are one week after a change in preventive measures as indicated in Table~\ref{tab_timeline}: July~8, August~13, August~25, October~6, October~21, November~11, December~21, and December~31 (see the Supplementary Materials for the results of all the other dates). For each of these eight dates, Figure~\ref{fig_scatter} shows for all 355 Dutch municipalities, each of which represented by one dot, the forecast and actual total number of reported positive tests within that municipality. Moreover, Figure~\ref{fig_maps} shows the forecast error and, as a scale reference, the total number of reported positive tests for these dates.

These figures confirm the observation made earlier on the national level that our forecasts are relatively accurate for periods where the transmission rates are stable. For instance, this is the case for the forecast for November 11-25 but not for the forecast for October~21 - November~4. Furthermore, the maps suggest that, in most cases, there are only a few municipalities with a relatively large forecast error. Thus, generally, the forecast error is distributed quite evenly over the municipalities. A notable exception to this is the forecast for October~21 - November~4, which can be explained by the relative instability of the corresponding transmission rates.

The graphs in Figure~\ref{fig_scatter} suggest that the variation in forecast error is quite high for July~8, August~13, and August~25, and relatively small for the other five dates. One explanation for this is that the actual numbers of daily new infections was relatively small compared to the other dates. In particular, many municipalities had only a few or even zero reported cases within the considered time periods, meaning that adjusting a small forecast error for population size may result in large differences per 100,000 inhabitants.

The graphs in Figure~\ref{fig_maps} do not suggest a strong relation between a large number of infections and a large forecast error. For instance, for October~21, large forecast errors occur for ranges of municipalities with both low incidence numbers (e.g., in the northern municipalities) or high incidence numbers (e.g., in the central eastern municipalities).

\begin{figure}[ht!]
    \centering
    \includegraphics[scale = .7]{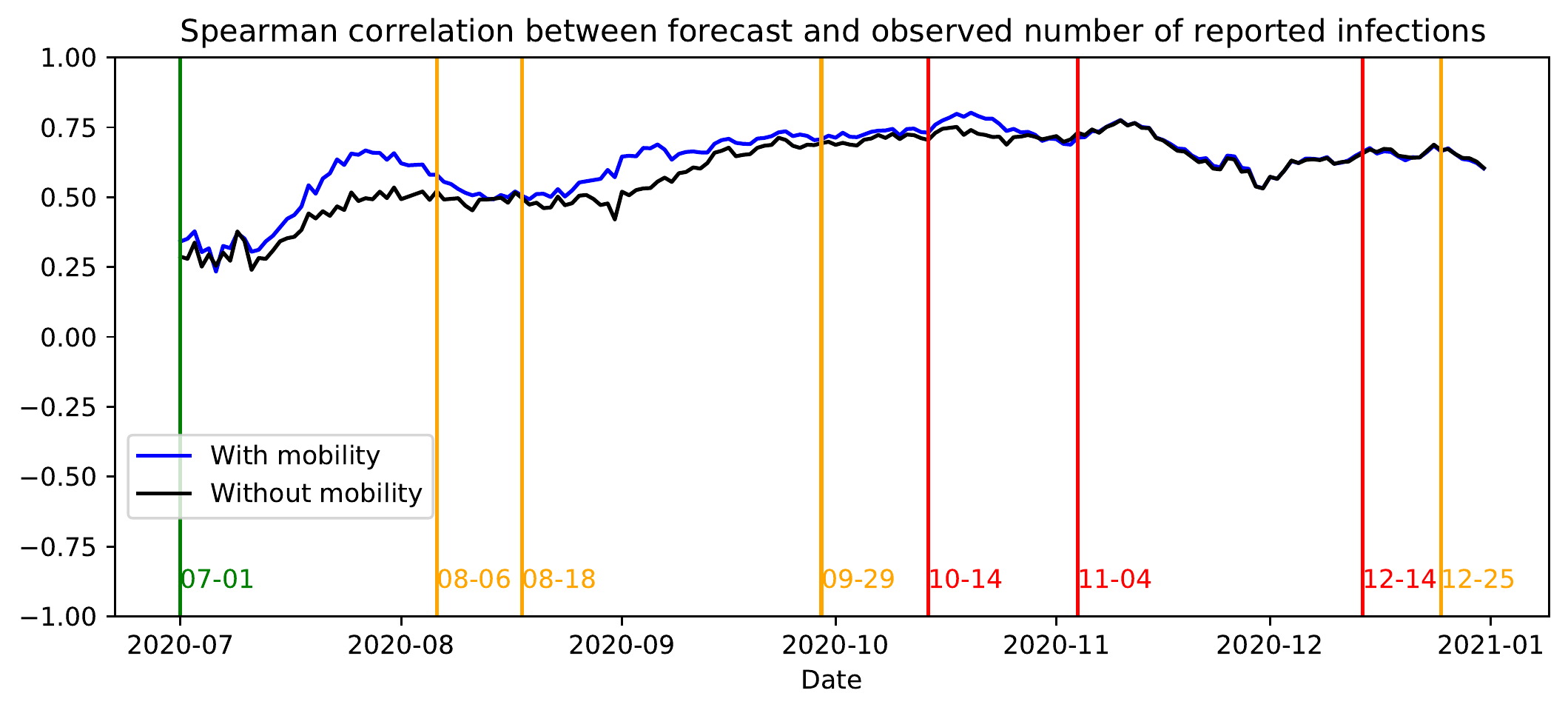}
    \caption{Correlation between the orders of the 14-day ahead forecast and observed fractions of reported infections.}
    \label{fig_final_frac_corr_14}
\end{figure}

We now assess the quality of our relative forecasts, i.e., forecasts of the fraction of reported infections in a given municipality compared to the nation-wide number. The blue line in Figure~\ref{fig_final_frac_corr_14} shows the Spearman correlation between the forecast and actual numbers of reported infections. Note that a correlation of 1 or $-1$ means that the forecast fractions have the same or reverse order as the observed fractions, respectively. The calculated values in Figure~\ref{fig_final_frac_corr_14} indicate that, apart from a small period at the start of July, the two orders are moderately and sometimes even strongly correlated (from mid-September until mid-November). This suggests that our approach is able to forecast the order of municipalities in terms of new infections with reasonable accuracy.

\subsection{Influence of mobility}
 
\begin{figure}[ht!]
\centering
\includegraphics[scale = .7]{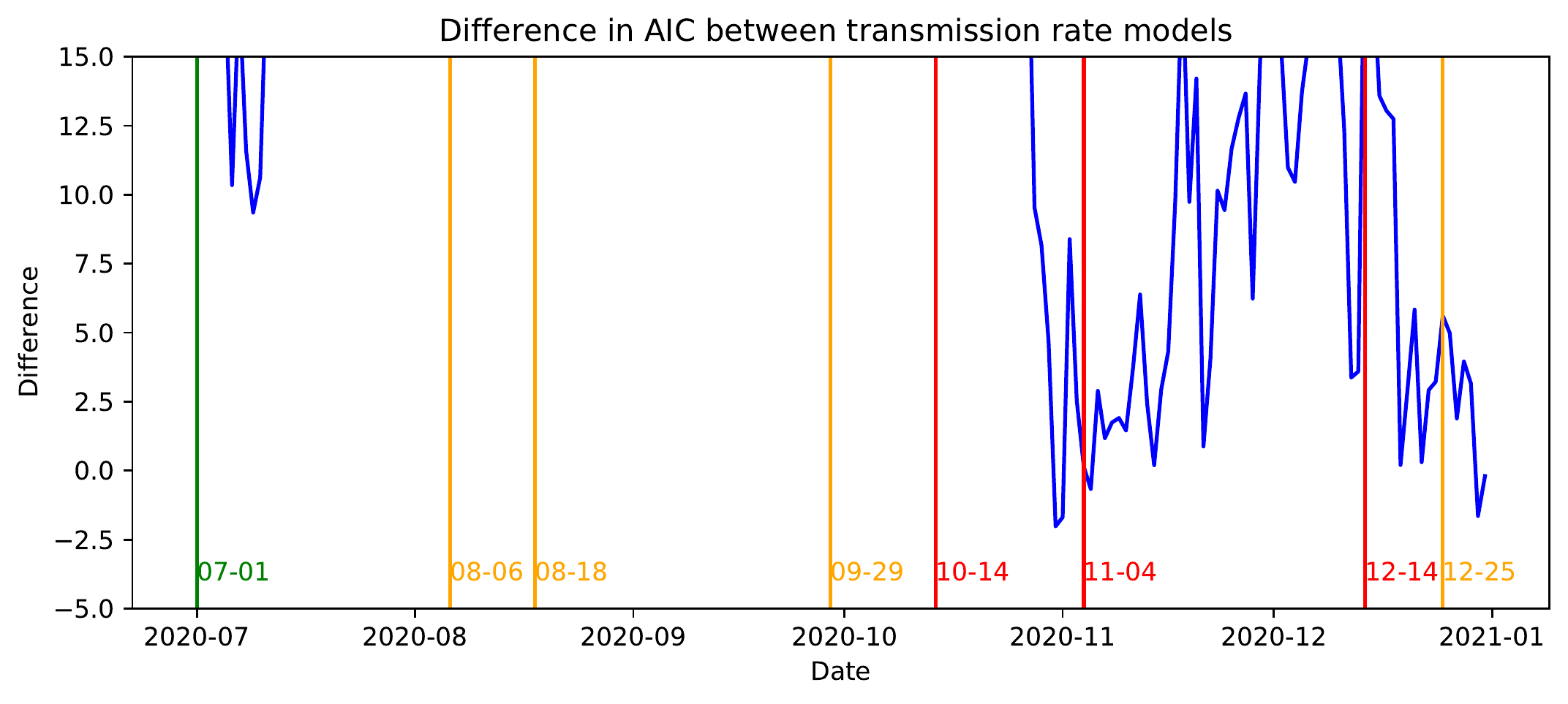}
\caption{The difference in AIC for the negative binomial models with and without mobility (zoomed in for differences between $-5$ and $15$).}
\label{fig_final_zoom_AIC}
\end{figure}

Finally, we focus on the difference in model performance and forecast accuracy between the models with and without mobility. First, Figure~\ref{fig_final_zoom_AIC} shows the differences in AIC between the negative binomial model with and without mobility. For clarity, we have zoomed in to a small band around the critical zero point that determines which of the models has a better fit. These results show that the differences in AIC between the negative binomial models with and without mobility are generally larger than $10$ before the partial lockdown initiated on October~14. This suggests that, until that day, incorporating mobility in the estimation process leads to a significantly better fit. However, throughout the remainder of the year, the differences are generally between $-2.5$ and $10$, meaning that for this period the fit may not significantly improve when mobility is taken into account. Note that this period matches almost perfectly with the time when there was an urgent advice to work from home. Therefore, these observations suggest that taking into account information on commuting mobility does not lead to a significantly better model fit when there is such a strong advice in place.

\begin{figure}[ht!]
    \centering
    \includegraphics[scale = .7]{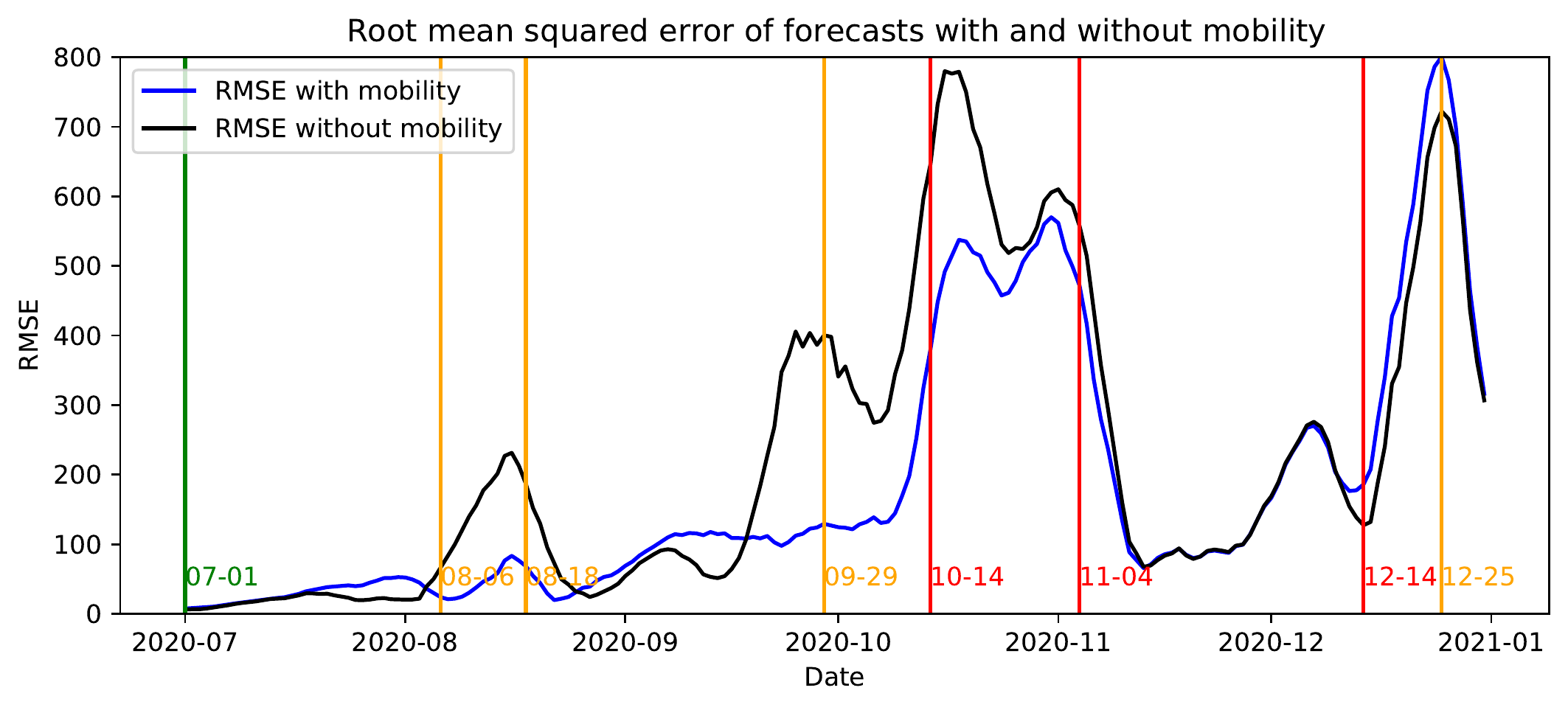}
    \caption{RMSE of the 14-day ahead forecasts (not per 100,000) of the models with and without mobility.}
    \label{fig_final_RMSE}
\end{figure}

To assess the additional value of taking mobility into account when forecasting the absolute numbers of new infections, we calculate for each date the RMSE of the corresponding forecast for both models and compare these values in Figure~\ref{fig_final_RMSE}. These results indicate that in most cases the forecasts that incorporate mobility are more accurate than those that do not. Moreover, after the partial lockdown measures of November 4, the RMSE of both models is practically equal. This suggests that the influence of mobility on forecast accuracy is insignificant in this period. 

With regard to the forecast accuracy in the relative number of new infections, we again consider the Spearman correlations in Figure~\ref{fig_final_frac_corr_14}. In essentially all cases, the correlation for the model with mobility is higher than that for the model without mobility. This means that taking mobility into account generally leads to better forecasts of the relative order of municipalities in terms of new infections. This suggests that taking mobility into account also improves the ability of our approach to identify local outbreaks. Furthermore, the correlations are practically the same in the period November-December, meaning that taking mobility into account does not significantly improve the forecasts in this period. This is consistent with the behavior of the RMSE as discussed above.

\section{Conclusions and discussion}
\label{sec_conclusion}

In this paper, we studied the transmission of SARS-CoV-2 in the Netherlands during the first year of the pandemic and the role of mobility on this transmission. To model transmissions between different regions, we employed the mobility-enhanced SEIR-compartmental model in \cite{Gosgens2021}. We obtained suitable parameter values for this model using commuter data and information on the seroprevalence of SARS-CoV-2. In particular, we developed a maximum likelihood estimation approach to fit suitable transmission rates for the compartmental model. Using the initialized model, we are able to provide accurate forecasts of the development of the transmission on the level of municipalities for up to 14 days. We found that taking mobility into account generally improves the model fit and forecast quality, but that this improvement is insignificant when mobility and contacts are restricted by national preventive measures. More precisely, we observed that the partial lockdown measures of November 4, 2020 led to sudden changes in parameter estimates and forecast accuracy. These observations suggest that these measures were particularly effective in reducing contacts and more so than those of the first partial lockdown on October 14, 2020. In particular, the closing of public places might have forced a significant part of the working population to work from home again. Such a change in mobility would be detected directly by the current model since the used mobility data consists purely of commuting information. 

We observe that the estimated parameters and forecast accuracy of the model generally remain relatively stable during periods where no new measures are being enforced. However, our work also shows that this changes dramatically when there is a change in mobility-restricting measures since the historic data used for initialization is no longer representative for the future. One way to correct for this, which we also aim to address in future work, is to integrate the most recent information on people's behavior into the estimation procedure, which could be obtained from, e.g., contact tracing apps.

We conclude by discussing several limitations of the research in this paper. First, we focus solely on inter-regional mobility and do not take other traits into account that influence susceptibility, infectivity, contact patterns and virus transmission, such as age, behavior, or household composition, into account. Also, we focus on the early phases of the pandemic, when vaccines were not yet available and hence the proportion of people who were protected was relatively small and only resulting from natural infection. It would be easy to add vaccination to the model, but it would add yet another layer of complexity because one would need to include how vaccinated individuals change their susceptibility, infectivity, contacts and mobility. Additionally, pandemic and as vaccinations increased, adherence to preventive measures changed in the population, as well as testing strategy and willingness to get tested.

Second, we use only commuter travel information as a mobility proxy. This is a limitation because we cannot test hypotheses and draw conclusions on general relations between mobility and, say, forecast quality. We also do not take other reasons for travel and travel destinations other than workplaces into account where transmission can occur such as social gatherings, mass events and holidays. To obtain a more realistic proxy for the mobility data, sources other than commuter information should be incorporated as well.

One option for this is travel behavior obtained via track-and-trace apps or location data obtained from cellphones \citep{Grantz2020}. However, these data typically contain privacy-sensitive information, for which the risk exists that observed movements are traced back to individual users. Therefore, these data must be sufficiently anonymized and/or aggregated to circumvent this issue. However, even when the data is anonymized, often permission for its usage is required. Summarizing, it is difficult, if not impossible, to obtain such detailed data. Therefore, one reasonable question for future research is how alternative data sources may be used to obtain suitable proxies for the actual mobility.

Finally, a similar problem occurs when certain epidemiological data is not available or no longer representative. For instance, an increase in vaccinated individuals leads to a decrease in hospitalizations. As a consequence, the difference between the number of new infections and new hospitalizations becomes larger, meaning that the latter becomes less representative of the number of infectious people. Since the latter is used to initialize a crucial parameter in our model (namely the fraction of positively tested infectious people), we expect the performance of the model and our initialization procedure in its current form to decrease when applied to a time period when a substantial part of the population has been vaccinated. One way to solve this issue could be to estimate the number of infectious people using other statistics that are known to correlate with this information, such as the daily fraction of positive tests.

\noindent \section*{Acknowledgments}

We thank Jetty Komrij for her input on the information needs of policy makers, Richard Post for his help with formulating the parameter estimation problem, and Nic Saadah for sharing his insights on human behavior during the COVID-19 pandemic. This work is supported by Netherlands Organisation for Scientific Research (NWO) through ZonMw grant no. 10430032010011. M.G., R.v.d.H., and N.L. are also supported by NWO through Gravitation NETWORKS grant no. 024.002.003.

\bibliographystyle{vancouver}
\bibliography{references}
\end{document}